\documentclass[10pt,reqno]{amsart}

\usepackage{amssymb,amsmath,amsthm,latexsym,color}
\usepackage{pifont}
\usepackage{verbatim}
\usepackage{amsfonts}
\usepackage{hyperref}
\usepackage[title,titletoc,toc]{appendix}

\title[Blow-up for the Hardy-H\'enon Equation]{Blow-Up and Life-Span for  Nonlinear Hardy-Hénon Parabolic Equations}
\author[B. Ben Slimene and S. Tayachi] {Byrame Ben Slimene and Slim Tayachi}
\address{%\vspace{0cm}\newline 
Byrame Ben Slimene, %\newline 
Taibah University, College of Science, Department of Mathematics, 42353 AL-Madinah, Saudi Arabia, 
e-mail: {\tt byramebenslimene@yahoo.fr}} %\vspace{1cm}
\address{Slim Tayachi, %\vspace{0cm}\newline 
Universit\'e de Tunis El Manar, Facult\'e des Sciences de Tunis, D\'epartement de
Math\'ematiques, Laboratoire  \'Equations aux D\'eriv\'ees Partielles LR03ES04,  2092 Tunis, Tunisia, 
e-mail: {\tt slim.tayachi@fst.rnu.tn} }
\subjclass[2010]{Primary 35K55, 35B44; Secondary 35K57, 35C15.} \keywords{Hardy-Hénon parabolic equations, finite time blow-up, Life-span.} 
\font\TenEns=msbm10 \font\SevenEns=msbm7 \font\FiveEns=msbm5
\newfam\Ensfam

\textfont\Ensfam=\TenEns \scriptfont\Ensfam=\SevenEns
\scriptscriptfont\Ensfam=\FiveEns
\def\N {\mathbb N}
\def\R {\mathbb R}

\def\Rd {{\mathbb R}^N}

\def\X {{\mathcal X}}

\newtheorem{To}{Theorem}[section]
\newtheorem{prop}[To]{Proposition}

\newtheorem{lem}[To]{Lemma}
\newtheorem{rem}{Remark}

\date{\today}
\begin{document}
\maketitle
\begin{abstract}
We consider the nonlinear Hardy-Hénon parabolic equation $
u_t=\Delta u+|x|^l|u|^\alpha u,$ $t>0,$ $x\in\mathbb{R}^N$ or $x\in\mathbb{R}^N\setminus\{0\}$,
with $N\ge1$, $l>-\min\{2,N\}$, $\alpha>0$, with nonnegative initial datum $u_0$. We first establish new criteria for finite-time blow-up. We then derive optimal upper bounds for the life-span $T_{\max}(\lambda\varphi)$ of solutions arising from the initial datum $u_0=\lambda\varphi$. These estimates establish a precise dependence of $T_{\max}(\lambda\varphi)$ on both the initial profile $\varphi$ and the parameter $\lambda$ in the small and large-amplitude regimes.
\end{abstract}

\medskip

\section{Introduction and main results}
\setcounter{equation}{0}
In this paper, we consider the semilinear Hardy-Hénon parabolic equation
\begin{equation}\label{NLHH}
u_t = \Delta u +|x|^l |u|^\alpha u,
\end{equation}
with initial data
\begin{equation}\label{initial data}
u(0)=u_0 \geq 0,
\end{equation}
where $u=u(t,x)$, $(t,x)\in (0, T)\times D,$  $T>0$, $D=\mathbb{R}^N$ if $l\ge0$, $D=\mathbb{R}^N\setminus\{0\}$ if $l<0,$ $N\in \N$,  $l>-\min\{2,N\}$ and    $\alpha>0$. Equation (\ref{NLHH}) with $l<0$ is known in the literature as the  Hardy parabolic equation. For $l>0$ it is known as the H\'{e}non parabolic equation. The case $l=0$ corresponds to the standard nonlinear heat equation. The problem \eqref{NLHH}--\eqref{initial data} is locally well posed in various Banach spaces; see, for example, \cite{BenSlimene, BTW, CFL, Ch, CIT2, CIT, CITT, CITT2, Tayachi, TW4, Wang} and the references therein. We denote by $T_{\max}(u_0)$ the maximal existence time of the solution to \eqref{NLHH}--\eqref{initial data}. If $0<T_{\max}(u_0)<\infty$, the solution is said to blow-up in finite time, whereas if $T_{\max}(u_0)=\infty$, it is called a global solution. By convention, we set $T_{\max}(u_0)=0$ whenever there is no local solution. The life-span of a solution refers to the maximal existence time corresponding to initial data of the form $u_0=\lambda\varphi$, where $\varphi$ is a fixed element of the Banach space under consideration, and $\lambda>0$.

For the case $l=0$, the life-span of solutions has been extensively studied; see, for example, \cite{GW,LeeNi,MY,TW4,TW5}. In the case $l\neq0$, lower bounds for the life-span were established in \cite{TW4}. See also \cite{CFL,TW4} for upper bound estimates of the life-span for certain classes of initial data.

In this paper, we investigate finite-time blow-up solutions for the case $l\neq0$. More precisely, we establish blow-up criteria for \eqref{NLHH}--\eqref{initial data} and derive optimal upper bounds for the life-span of solutions corresponding to a broad class of nonnegative initial data $u_0=\lambda\varphi$, expressed in terms of the parameter $\lambda$ and $\varphi$.

As is standard practice, we study \eqref{NLHH}--\eqref{initial data} via the associated integral equation
\begin{equation}\label{NLHHint}
u(t) = e^{t\Delta} u_0 + \int_0^t e^{(t-\sigma)\Delta} \big[|\cdot|^l |u(\sigma)|^\alpha u(\sigma)\big]  d\sigma,
\end{equation}
where, for
all $t>0,$ $e^{t\Delta}$ denotes the heat semigroup
\begin{equation}
\label{semigroup1}
e^{t\Delta}f=G(t,\cdot)\ast f,
\end{equation}
with $\ast$ denotes the convolution operator and
\begin{equation}
\label{semigroup2}
G (t,x)=(4\pi t)^{-\frac{N}{2}}e^{-{|x|^2\over 4t}},\; t>0,\; x\in \R^N.
\end{equation}
A solution $u$  of \eqref{NLHHint} is a
mild solution of the problem \eqref{NLHH}--\eqref{initial data}.

In what follows, for $1 \leq q \leq \infty$, we denote by $L^q(\mathbb{R}^N)$ the usual Lebesgue space endowed with its standard norm $\|\cdot\|_{L^q}$. We also consider the weighted Lebesgue space $L^q_\gamma(\mathbb{R}^N)$, consisting of all measurable functions $f : \mathbb{R}^N \to \mathbb{R}$ such that $|\cdot|^\gamma f \in L^q(\mathbb{R}^N)$, and equipped with the norm $\|f\|_{L^q_\gamma} = \||\cdot|^\gamma f\|_{L^q}$, where $\gamma \in \mathbb R$. Finally, $C_0(\mathbb{R}^N)$ denotes the space of continuous functions vanishing at infinity, while $C_b(\mathbb{R}^N)$ denotes the space of continuous  bounded functions on $\mathbb{R}^N$. 

By Fujita-type results \cite{BL, F, Fujita2, LM, Qi}, when $\alpha \leq \frac{2+l}{N}$,  problem (\ref{NLHH})-(\ref{initial data})  admits no nontrivial nonnegative global solution for any nontrivial  nonnegative initial data, while if $\alpha > \frac{2+l}{N}$, there are global positive solutions. In the supercritical case $\alpha > \frac{2+l}{N}$, additional conditions on the initial data $u_0$ may still lead to a finite-time explosion. To illustrate this phenomenon, first consider the case $l=0$.
It is proved in \cite[Theorem 3.2 (i), p. 372]{LeeNi} that for $\alpha > \frac{2}{N}$, there exists $C>0$ such that if  $u_0\in C_b(\R^N),$ $u_0\geq 0$ and  $\liminf_{|x|\to \infty} |x|^{\frac{2}{\alpha}} u_0(x) > C$, then $T_{\max}(u_0) < \infty$. For instance,  for all $\alpha>0,$ if $u_0\in C_b(\R^N),$  $u_0\geq 0$ and $\liminf_{|x|\to \infty} |x|^{\gamma}u_0(x)>0$ with $\gamma<{2\over\alpha},$ then $T_{\max}(u_0)<\infty,$ as shown in \cite{LeeNi}. These results have been significantly refined in many works; see \cite{ MY,  R, SW, TW5}. In particular, the authors of \cite{TW5} extended the analysis to a broader class of initial data not covered by previous results. For instance, they showed that if $u_0 \in C_0(\mathbb{R}^N)$, $u_0 \geq 0$, and $u_0(x) \geq c |x|^{-2/\alpha} \sin^2(\log |x|)$ for all sufficiently large $|x|$, with $c>0$ large enough, then $T_{\max}(u_0) < \infty$. Fewer results are available in the case $l \neq 0$. In the Hardy case, that is $-\min\{2,N\} < l < 0$, it is observed in \cite[Remark 1.3, p. 119]{BTW} that if $N \geq 3$ and $\frac{2+l}{N} < \alpha < 2 \frac{2+l}{N-2}$, then there exists $C>0$ such that if $u_0 \in C_0(\mathbb{R}^N)$, $u_0>0$, and $\liminf_{|x| \to \infty} |x|^{\frac{2+l}{\alpha}} u_0(x) \geq C$, one has $T_{\max}(u_0) < \infty$.  On the other hand, when $l \geq 0$, corresponding to the Hénon case, Wang \cite[Theorem 0.2', p. 556]{Wang} proved that for $N \geq 3$, if $u_0 \in C_b(\mathbb{R}^N)$ satisfies $0 \leq u_0(x) \leq c |x|^{-\frac{l}{\alpha}}$ in $\mathbb{R}^N$, $\alpha \ge  \alpha_c(l)$, (where $\alpha_c(l)$ denotes the Joseph-Lundgren exponent defined in \cite[equation (0.7), p. 555]{Wang},) and $\liminf_{|x| \to \infty} |x|^{\frac{2+l}{\alpha}} u_0(x)$ is sufficiently large, then the corresponding solution of \eqref{NLHH}--\eqref{initial data} blows up in finite time.  It was recently proved in  \cite[Theorem 5.1 (ii)]{CFL} that, if $l>-\min \{2, N\}$, $\alpha>\frac{2+l}{N}$ and $0<\gamma<\frac{2+l}{\alpha}$, then, for any $u_0 \in C_0(\mathbb R^N)$ satisfying $u_0 \ge 0$ and $\liminf_{|x|\to\infty}|x|^\gamma u_0(x)>0$, we have $T_{\max}(u_0) < \infty$.

Throughout this paper, we let $D_\lambda$ denote the dilation operator
\begin{equation}
\label{dilation}
(D_\lambda f)(x) = f(\lambda x), \quad \forall\; \lambda > 0.
\end{equation}
For $-\min\{2,N\} < l < N\alpha$, we define
\begin{equation}\label{constanttype}
C_l = (4\pi)^{- \frac{N}{2}} \int_{\mathbb{R}^N} e^{-\frac{|y|^2}{4}} |y|^{-\frac{l}{\alpha}}  dy= \left(e^{\Delta} |\cdot|^{-\frac{l}{\alpha}} \right)(0)=\frac{2^{-\frac{l}{\alpha}}}{\Gamma\left(\frac{N}{2}\right)}
\Gamma\left(\frac{N}{2}-\frac{l}{2\alpha}\right),
\end{equation}
where $\Gamma$ is the Euler gamma function. This constant will appear repeatedly in our analysis. 
\subsection{Main results on blow-up criteria}
Our main results on blow-up criteria, which extend from \cite[Theorem 1.7, p. 146]{TW5} to $l\not=0$, are as follows. The first theorem is for the case $l>0.$ Note that the blow-up may hold at the origin, as it may also not hold at the origin. See \cite{GSouplet,GLS1,GLS2,GS,FT}. 
\begin{To}[Blow-up criterion for the Hénon equation]\label{thmblowupHénon}
Let $N \geq 1$, and let $l,$ $\alpha,$ $\gamma$ be real numbers satisfying
\begin{equation*}
 0 < l <N\alpha,\quad {l \over \alpha}\leq \gamma<N, \quad 0 <\alpha \le {2+l\over \gamma}.
\end{equation*}
Let $w(x)=(1+ |x|)^{\frac{l}{\alpha}}$ and define the spaces $C_w(\mathbb R^N)=\left\{ \phi \in C(\mathbb R^N); \quad  \|\phi w\|_{L^\infty}<\infty \right\}$ and  $ L^\infty_{w}\left(\R^N\right)= \left\{ \phi; \quad  \|\phi w\|_{L^\infty}<\infty \right\}$. Given $u_0 \in C_w(\mathbb R^N)$, let $u\in L^\infty_{loc}\left(\left[0,T_{\max}(u_0)\right), L^\infty_{w}\left(\R^N\right)\right)$ denote the maximal solution of \eqref{NLHHint} (by \cite[Theorem 2.5, p. 565]{Wang}, this solution is in fact a classical solution of \eqref{NLHH}--\eqref{initial data}.) Let  $f\in C_w(\mathbb R^N)$ be such that
\begin{itemize}
\item[(i)] $0 \le f(x) \le u_0(x)$ for a.e. $x \in \mathbb{R}^N$,
\item[(ii)] there exists $K>0$ such that
\begin{equation} \label{initialinfmodulef} 
f(x) \le K |x|^{-\gamma} \quad \text{for a.e. } x \in \mathbb{R}^N.
\end{equation}
\end{itemize}
Suppose there exist a sequence $\lambda_n \to \infty$ and a function  $z \in L^1_{loc}\left({\mathbb R^N} \right)$, with $z \not\equiv 0$ a.e.,  such that
\begin{equation}\label{subseqlim}
\lambda_n^{\gamma} D_{\lambda_n} f \to z \quad \text{in } \mathcal{D}'(\mathbb{R}^N), \text{ as } n \to \infty.
\end{equation}
If $\alpha = \frac{2+l}{\gamma}$, assume in addition that
\begin{equation}\label{notsmallHenon}
\left\|e^{\Delta} z\right\|_{L^\infty}  > C_l \left( \frac{2+l}{2 \alpha} \right)^{\frac{1}{\alpha}},
\end{equation}
where $C_l$ is given by \eqref{constanttype}.  Then $T_{\max}(u_0) < \infty$.
\end{To}
\begin{rem} {\rm Some remarks are in order concerning Theorem \ref{thmblowupHénon}.
\begin{itemize}
\item[(i)]The condition ${l \over \alpha}\leq \gamma$ is used in \cite{TW4,CITT,CITT2} to construct local solutions. It seems also it is necessary for the existence of a solution. See \cite[Theorem 2.5 (ii), p. 565]{Wang}. See also \cite{IL} and references therein.
\item[(ii)] Conditions $\gamma\leq {2+l\over \alpha}$ and \eqref{notsmallHenon} are optimal. In fact for $|u_0(x)|\leq C{1\over (1+|x|^2)^{{\sigma\over 2}}}, \; {2+l\over \alpha} \leq \sigma<N,$ we have global existence, for $C>0$ sufficiently small, see \cite[Theorem 1.5, Remark 1.6 (a)]{CITT}.
\item[(iii)] The previous theorem holds for $l=0$ by supposing moreover $\gamma>0$.
\item[(iv)] To emphasize the dependence on $l$, we denote $T_{\max}(u_0)$ by $T_{\max}(l,u_0)$. Fix $\alpha>\frac{2}{N}$,  $ 0 < l <N\alpha$,  and let $\max\{0,l-2\}\le l_1<\min\{N\alpha-2,l\}$. Choose $u_0\in C_w(\mathbb R^N)$ such that $u_0(x)=c|x|^{-\frac{2+l_1}{\alpha}}$, for all $|x|\geq\rho$.  Taking  $f=u_0$, we have, for all $|x| \ge \frac{\rho}{\lambda}$, $\lambda^{\frac{2+l_1}{\alpha}} D_\lambda f(x) = c \lambda^{\frac{2+l_1}{\alpha}} |\lambda x|^{-\frac{2+l_1}{\alpha}} =  c |x|^{-\frac{2+l_1}{\alpha}}$. Consequently,  $\lambda^{\frac{2+l_1}{\alpha}} D_\lambda f \to c |\cdot|^{-\frac{2+l_1}{\alpha}}$ in $\mathcal{D}'(\mathbb{R}^N)$ as $\lambda \to \infty$. Therefore, by Theorem \ref{thmblowupHénon}, we obtain $T_{\max}(l,u_0)<\infty$. On the other hand, since $\frac{2}{\alpha} \le\frac{2+l_1}{\alpha}<N$, it follows from  \cite[Theorem 1.5, Remark 1.6 (a)]{CITT} that $T_{\max}(0,u_0)=\infty$, provided that the additional condition  $|u_0(x)| \le 2 c \frac{1}{(1+ |x|^2)^{\frac{2+l_1}{2 \alpha}}}$ is satisfied, where $c>0$ is sufficiently small. This example shows that the same initial datum can generate a global solution for the standard nonlinear heat equation, while the corresponding equation with the Hénon weight admits a solution that blows up in finite time. 
\end{itemize}}
\end{rem}
The second theorem is for the case $l<0.$
\begin{To}[Blow-up criterion for the Hardy equation]\label{thmblowup}
Let $N \geq 1$, and let $l,$ $\alpha,$ $\gamma$ be real numbers satisfying
\begin{equation*}
   -\min\{2,N\}<l <0,\quad 0<\gamma<N, \quad 0 <\alpha \le {2+l\over \gamma}.
\end{equation*}
Let $u_0 \in C_0(\mathbb{R}^N)$, and let $u \in C\left([0, T_{\max}(u_0)), C_0(\mathbb{R}^N)\right)$ be the maximal solution of \eqref{NLHHint}. Let $f \in C_0(\mathbb{R}^N)$ be a function such that 
\begin{itemize}
\item[(i)] $0 \le f(x) \le u_0(x)$ for a.e. $x \in \mathbb{R}^N$,
\item[(ii)] there exists $K>0$ such that
\begin{equation} \label{initialinfmodulefz} 
f(x) \le K |x|^{-\gamma} \quad \text{for a.e. } x \in \mathbb{R}^N.
\end{equation}
\end{itemize}
Assume there exist a sequence $\lambda_n \to \infty$ and a function  $z \in L^1_{loc}\left({\mathbb R^N} \right)$, with $z \not\equiv 0$ a.e.,  such that \begin{equation}\label{subseqlimz}
\lambda_n^{\gamma} D_{\lambda_n} f \to z \quad \text{in } \mathcal{D}'(\mathbb{R}^N), \text{ as } n \to \infty.
\end{equation}
If $\alpha = \frac{2+l}{\gamma}$, we further suppose that
\begin{equation}\label{notsmall}
\left(e^{\Delta} z\right)(0) > C_l \left( \frac{2+l}{2 \alpha} \right)^{\frac{1}{\alpha}},
\end{equation}
where $C_l$ is defined by \eqref{constanttype}.  Then $T_{\max}(u_0) < \infty$.
\end{To}

\begin{rem}\label{compact}$\;${\rm
\begin{itemize}
\item[(i)] Although condition $\eqref{notsmall}$ implies that $\|e^\Delta z\|_{L^\infty}> C_l\left(\frac{2+l}{2\alpha}\right)^{\frac{1}{\alpha}}$, we impose the stronger condition $\eqref{notsmall}$ when $l<0$, rather than requiring only the above inequality, as in the case $l\geq 0$. One possible interpretation of this stronger assumption is that it ensures that the required largeness holds at the origin, where the potential $|x|^l$  is singular.
\item[(ii)] Conditions $\gamma\leq {2+l\over \alpha}$ and \eqref{notsmall} are optimal. In fact for $|u_0(x)|\leq C{1\over (1+|x|^2)^{{\sigma\over 2}}}, \; {2+l\over \alpha} \leq \sigma<N,$ we have global existence, for $C>0$ sufficiently small, see \cite[Theorem 1.3, p. 119]{BTW}.
\item[(iii)] Fix $\alpha>\frac{2+ l}{N}$, $-\min\{2,N\}<l <0$, and let $l\le l_2< \min\{N \alpha-2, 0\}$. Choose $u_0\in C_0(\mathbb R^N)$ such that $u_0(x)=c|x|^{-\frac{2+l_2}{\alpha}}$, for all $|x|\geq\rho$. Let $f=u_0$. Then, for all $|x|\geq \frac{\rho}{\lambda}$, we have $\lambda^{\frac{2+l_2}{\alpha}}D_\lambda f(x)=c|x|^{-\frac{2+l_2}{\alpha}}$. Consequently, $\lambda^{\frac{2+l_2}{\alpha}}D_\lambda f\to c|\cdot|^{-\frac{2+l_2}{\alpha}}$ in $\mathcal D'(\mathbb R^N)$ as $\lambda\to\infty$. It follows from \cite[Theorem 1.7, p. 146]{TW5} that $T_{\max}(0,u_0)<\infty$. On the other hand, by Part (ii) of this remark, we have $T_{\max}(l,u_0)=\infty$, provided that  $|u_0(x)| \le  \frac{2  c}{(1+ |x|^2)^{\frac{2+l_2}{2 \alpha}}}$, where $c>0$ is sufficiently small. Thus, the same initial datum can lead to finite-time blow-up for the equation without the Hardy weight, while giving rise to a global solution when the Hardy weight is present.  
\item [(iv)]  In the proofs of Theorems \ref{thmblowupHénon} and \ref{thmblowup}, it suffices to take $u_0 = f$. Both results extend to any $u_0 \in L^1_{\mathrm{loc}}(\mathbb{R}^N)$ satisfying $u_0 \ge f$ almost everywhere, provided that a local solution to \eqref{NLHHint} exists. For example, we may take $f\in L^\infty_{\gamma}(\R^N)$ and ${l\over \alpha}\leq \gamma<{2+l\over \alpha},\; 0<\gamma<N.$ See \cite{CITT,CITT2,TW4}.
 \item [(v)] It follows from  \cite[Proposition 2.1, p. 1110]{CDW-DCDS} that for any sequence $\lambda_n \to \infty$ there exists a subsequence (still denoted $\lambda_n$)  and a function $z \in L^1_{loc}\left(\mathbb{R}^N \right)$ such that \eqref{subseqlim} holds.  The main issue in Theorems \ref{thmblowupHénon} and  \ref{thmblowup} is that we require $z \not\equiv 0$,  or, if $\alpha = \frac{2+l}{\gamma}$, that $z$ is sufficiently large in a suitable sense. 
     \item [(vi)] Suppose that  the hypotheses of Theorem \ref{thmblowupHénon} or Theorem \ref{thmblowup} are satisfied and that $\alpha<\frac{2+l}{\gamma}$. Let us introduce the set
   \begin{equation*}
     I_\gamma=\{\phi\in C_0(\mathbb R^N); \quad \liminf_{|x|\to\infty}|x|^\gamma\phi(x)>0\}.
   \end{equation*}
    It is easy to see that if $u_0\in I_\gamma$, then there exist constants $C>0$ and $\rho>0$ such that $u_0(x)\ge C|x|^{-\gamma}$, for all  $|x| \ge \rho$. Taking $f\in C_0(\mathbb R^N)$ such that $0\leq f\leq u_0,$ $f(x)= C|x|^{-\gamma}$ for $|x|\ge  \rho$.   We have, for all $|x| \ge \frac{\rho}{\lambda}$, $\lambda^{\gamma} D_\lambda f(x) =  C |x|^{-\gamma}$, and thus $\lambda^{\gamma} D_\lambda f \to C |\cdot|^{-\gamma}$ in $\mathcal{D}'(\mathbb{R}^N)$ as $\lambda \to \infty$. Consequently, Theorem \ref{thmblowupHénon} or Theorem \ref{thmblowup}, as appropriate, implies that  every solution with nonnegative initial data  $u_0\in I_\gamma$ blows up in finite time. This recovers Theorem 5.1 (ii) of \cite{CFL}. Our results, however, apply to a strictly larger class of initial data. Indeed, even within $C_0(\mathbb R^N)$, there exist initial data satisfying the hypotheses of Theorem \ref{thmblowupHénon} or Theorem \ref{thmblowup} that do not belong to $I_\gamma$. To illustrate this point, consider a nonnegative function $u_0\in C_0(\mathbb R^N)$ such that, for some $\rho>0$, $u_0(x)=|x|^{-\gamma}\sin^2(\log|x|)$ for $|x|\ge\rho$. We have  $\liminf_{|x|\to\infty}|x|^\gamma u_0(x)=0$, and hence $u_0\notin I_\gamma$.  Moreover, if $f\in C_0(\mathbb R^N)$ such that $0\leq f\leq u_0,$ $f(x) =  |x|^{-\gamma}\sin^2(\log|x|)$ for $|x|\ge  \rho$, then $\lambda^{\gamma} D_\lambda f(x) = |x|^{-\gamma} \sin^2(\log |x| + \log \lambda)$ for all $|x| \ge \frac{\rho}{\lambda}$, and in this case we may choose $\lambda_n = e^{n \pi}$ to apply the preceding results. Therefore, this initial datum  $u_0$ generates a solution that blows up in finite time.
 \item [(vii)] Let $\alpha = \frac{2+l}{\gamma}$. Suppose $f \in  C_w(\mathbb R^N)$ or $f \in C_0(\mathbb{R}^N)$  is nonnegative and satisfies $f(x) = c |x|^{-\frac{2+l}{\alpha}}$ for all $|x| \ge \rho>0,$ and with $c>0$ is large enough. Then, for every $|x| \ge \frac{\rho}{\lambda}$, $\lambda^{\frac{2+l}{\alpha}} D_\lambda f(x)  = c |x|^{-\frac{2+l}{\alpha}}$, and thus $\lambda^{\frac{2+l}{\alpha}} D_\lambda f \to c|\cdot|^{-\frac{2+l}{\alpha}}$ in $\mathcal{D}'(\mathbb{R}^N)$ as $\lambda \to \infty$, so that the above results apply to any $u_0 \ge f$. Moreover, if $f(x) = c |x|^{-\frac{2+l}{\alpha}} \sin^2(\log |x|)$ for $|x| \ge \rho$, then $\lambda^{\frac{2+l}{\alpha}} D_\lambda f(x) = c |x|^{-\frac{2+l}{\alpha}} \sin^2(\log |x| + \log \lambda)$ for $|x| \ge \frac{\rho}{\lambda}$, and in this case choosing $\lambda_n = e^{n \pi}$ allows us to apply the previous results to any  $u_0 \ge f$. Note that the critical case $\alpha = \frac{2+l}{\gamma}$ is not considered in  \cite[Theorem 5.1 (ii)]{CFL}.
\end{itemize}}
\end{rem}

\subsection{Main results on the upper bounds of the life-span.}
We now turn to life-span results. Most of the existing results are related to the standard heat equation corresponding to $l = 0$. These results investigate the asymptotic behavior of $T_{\max}(\lambda \varphi)$ for nontrivial, nonnegative initial data $\varphi$ such as $\lambda \to \infty$ or $\lambda \to 0$. The behavior strongly depends on the class to which $\varphi$ belongs. For example, the asymptotics of $T_{\max}(\lambda \varphi)$ as $\lambda \to 0$ differ significantly depending on whether $\varphi \in C_b(\mathbb{R}^N)$ or $\varphi$ is singular at the origin. Results on the asymptotic behavior of the life-span for the heat equation can be found in \cite{GW, HisaIshige, LeeNi,  Tayachi1, TW4, TW5} and the references therein. For the general case $l> - \min\{2, N\}$, Pinsky in \cite{Pinsky} considered two classes of initial data $\varphi \in C_b(\mathbb R^N)$ with $\varphi \geq 0$. The first class consists of $\varphi$ dominated by a Gaussian, and the second consists of $\varphi$ such that there exist constants $c_1$, $c_2 > 0$ with $c_1 \le \varphi \le c_2$. He studied the asymptotic behavior of $T_{\max}(\lambda \varphi)$ as $\lambda \to 0$ when $T_{\max}(\lambda \varphi) < \infty$ for all $\lambda > 0$. The asymptotic order is shown to depend strongly on $l$, $\varphi$, $\alpha$ and $N$.  The authors of \cite{TW4} later established several lower bound estimates for the life-span for initial data belonging to Lebesgue and weighted Lebesgue spaces, as well as $C_0(\mathbb{R}^N)$, without assuming non-negativity of the initial data. Since these results are necessary for obtaining the asymptotic behavior of the life-span in this paper, we briefly recall the results needed for our study.
\begin{To}[Nonlinear Hardy parabolic equation {\cite[Theorem A.1]{TW4}}]\label{thm lower bound 1}
Let $N\ge 1$, $\alpha>0$ and $- \min \{2, N\}<l<0$. Let $\varphi \in  C_0(\mathbb R^N)$, and let $u \in C\left( [0, T_{\max}(\lambda \varphi)); C_0(\mathbb R^N)\right)$ be the maximal solution of \eqref{NLHHint} with initial data $\lambda \varphi$, $\lambda>0$. Then there exists a constant $C=C(N, l, \alpha)>0$ such that 
\begin{equation}\label{eqlower}
  T_{\max}(\lambda \varphi) \ge  C \left( \lambda \|\varphi\|_{L^\infty}\right)^{-\frac{2 \alpha}{2+ l}}, \quad \forall \lambda>0.
\end{equation}
\end{To}
\begin{To}[Nonlinear H\'enon parabolic equation {\cite[Theorem B.1 and Corollary B.2]{TW4}}]\label{thm lower bound 2bis}
Let $N\ge 1$, $\alpha>0$ and $0 \le l<N  \alpha$.  Let $\varphi \in L_{\frac{l}{\alpha}}^{\infty}(\mathbb{R}^N)$, and let $u \in C\left((0, T_{\max}(\lambda \varphi)); L_{\frac{l}{\alpha}}^{\infty}(\mathbb{R}^N)\right)$, $\lim_{t \to 0} \|u(t)- e^{t \Delta}\|_{L_{\frac{l}{\alpha}}^{\infty}}=0$  be the maximal solution of \eqref{NLHHint} with initial data $\lambda \varphi$, $\lambda>0$. Then there exists a constant $C = C(N, l, \alpha) > 0$ such that
\begin{equation}\label{eqlower2}
T_{\max}(\lambda \varphi) \ge C \left( \lambda \|\varphi\|_{L_{\frac{l}{\alpha}}^{\infty}}\right)^{-\alpha}, \quad \forall \lambda>0.
\end{equation}
\end{To}
\begin{To}[Nonlinear Hardy-Hénon parabolic equation {\cite[Theorem C.1 and Corollary C.2]{TW4}}]\label{thm lower bound 2}
Let $N\ge 1$, $\alpha>0$ and $- \min\{2, N\}<l<N  \alpha$. Let the real number $\gamma$ be such that $\max\{ 0, \frac{l}{\alpha}\} <\gamma< \min \{N, \frac{2+l}{\alpha} \}$. Let $\varphi \in L_{\gamma}^{\infty}(\mathbb{R}^N)$, and let $u \in C\left((0, T_{\max}(\lambda \varphi)); L_{\gamma}^{\infty}(\mathbb{R}^N)\right) \cap C\left((0, T_{\max}(\lambda \varphi)); L_{\nu}^{r}(\mathbb{R}^N)\right)$, $\lim_{t \to 0} \|u(t)- e^{t \Delta}\|_{L_{\gamma}^{\infty}}=0$ be the maximal solution of \eqref{NLHHint} with initial data $\lambda \varphi$, $\lambda>0$, where $r$ and $\nu$ are suitably chosen. Then there exists a constant $C = C(N, l, \alpha, \gamma) > 0$ such that
\begin{equation}\label{eqlower2bis}
T_{\max}(\lambda \varphi) \ge C \left( \lambda \|\varphi\|_{L_{\gamma}^{\infty}}\right)^{-(\frac{2+ l}{2 \alpha}-\frac{\gamma}{2})^{-1}}, \quad \forall\; \lambda>0.
\end{equation}
\end{To}

We now study upper bounds for the life-span. The upper bounds results for the life-span with $l\not=0$ are derived via scaling arguments in \cite[Proposition C.3, p. 621]{TW4}, but only under the restrictive assumptions  $l> -\min\{2, N\}$, $\max\{0, \frac{l}{N}\}<\alpha<\frac{2+l}{N}$ and $\max\{0, \frac{l}{\alpha}\}<\gamma< N$; moreover, the corresponding estimates do not provide explicit constants. Our results complete the study of the lower bounds for the life-span initiated in \cite{TW4} for (\ref{NLHH})–(\ref{initial data}), and allow us to analyze the asymptotic behavior of $T_{\max}(\lambda \varphi)$ both as $\lambda \to \infty$ and as $\lambda \to 0$. Moreover, as in \cite{TW4}, and in contrast to many life-span estimates in the literature, our upper bounds not only characterize the dependence on $\lambda$ but also provide a precise quantification of the dependence on the initial profile $\varphi$.

We distinguish between large and small values of $\lambda>0$. We begin with the case of large $\lambda$. The corresponding results are stated below.
\begin{To}\label{thmupper1bis}
Let $N \geq 1$, $\alpha > 0$ and $0 < l < N\alpha$. Let  $\varphi \in L_{\frac{l}{\alpha}}^\infty(\mathbb R^N)$, $\varphi \geq 0$ and $\varphi \not \equiv 0$. It follows that  $T_{\max}(\lambda \varphi) < \infty$ for all sufficiently large $\lambda > 0$, and 
\begin{equation*}
  \limsup_{\lambda \to \infty} \lambda^{\alpha} T_{\max}(\lambda \varphi) \le \frac{1}{\alpha} \left(\frac{ \| {e}^{\Delta} |\cdot|^{-\frac{l}{\alpha}} \|_{L_{\frac{l}{\alpha}}^\infty}}{ \|\varphi\|_{L_{\frac{l}{\alpha}}^\infty}}\right)^{\alpha}.
\end{equation*}
\end{To}
\begin{rem}{$\;${\rm
\begin{itemize}
  \item [1)] Let $N \geq 1$, $\alpha > 0$ and  $0 < l < N\alpha$.   Theorem \ref{thm lower bound 2bis} and Theorem \ref{thmupper1bis} together show that
\begin{equation}\label{Llinfty}
  T_{\max}(\lambda \varphi) \sim \lambda^{- \alpha}, \quad  \lambda \to \infty, 
\end{equation}
whenever $\varphi \in L_{\frac{l}{\alpha}}^\infty(\mathbb R^N)$, $\varphi \geq 0$ and $\varphi \not \equiv 0$.
  \item [2)] The equivalence \eqref{Llinfty} relies on the lower bound estimate \eqref{eqlower2}, which is obtained from the well-posedness theory in $L_{\frac{l}{\alpha}}^{\infty}(\mathbb{R}^N)$. This naturally raises the question of whether the problem \eqref{NLHHint} with $l>0$ is also well posed in the usual space $L^\infty(\mathbb{R}^N)$. The following argument shows that this is not the case. Let $N \geq 1$, $\alpha>0$ and $l>0$. Arguing by contradiction and assume that for every $u_0 \in L^\infty(\mathbb{R}^N)$ there exist $T(u_0)>0$ and a solution $u \in C\left((0,T(u_0)]; L^\infty(\mathbb{R}^N)\right)$  of \eqref{NLHHint}. Take $u_0 \equiv 1$. Then, for any $t \in (0,T(u_0)]$, \eqref{NLHHint} gives
      \begin{equation*}
        u(t) \ge \int_0^t e^{(t-\sigma)\Delta} |\cdot|^l  d\sigma.
      \end{equation*}
Using the heat kernel representation and the reverse triangle inequality, we obtain
\begin{equation*}
  u(t,x) \ge  (|x|-1)^l \int_0^t \int_{|z|\le 1} G(t-\sigma, z) dz d\sigma \to \infty, \quad \mbox{as} \quad |x|\to \infty,
\end{equation*}
which contradicts the fact that $u(t)\in L^\infty(\mathbb{R}^N)$.
\end{itemize}}}
\end{rem}
\begin{To}\label{thmupper1}
Let $N \geq 1$, $\alpha > 0$, $-\min\{2, N\}<l<0$, and let $C_l$ be given by \eqref{constanttype}. Let $\varphi \in C_0(\mathbb R^N)$, $\varphi \geq 0$ and $\varphi(0) \ne 0$. Then, for sufficiently large $\lambda > 0$, the  existence time of the maximal  solution of \eqref{NLHHint} with initial data $u_0 = \lambda \varphi$ satisfies  $T_{\max}(\lambda \varphi) < \infty$. In addition, one obtains the following upper bound:
\begin{equation}\label{ineq45}
  \limsup_{\lambda \to \infty} \lambda^{\frac{2 \alpha}{2+l}} T_{\max}(\lambda \varphi)  \leq  \left(\frac{2+l}{2 \alpha}\right)^{\frac{2}{2+ l}} \left(\frac{C_l}{\varphi(0)}\right)^{\frac{2 \alpha}{2+ l}}.
\end{equation}
\end{To}
\begin{rem}{\rm
Let $N \geq 1$, $\alpha > 0$ and $-\min\{2, N\}<l<0$. Combining Theorem \ref{thm lower bound 1} and Theorem \ref{thmupper1}, we obtain
\begin{equation*}
  T_{\max}(\lambda \varphi) \sim  \lambda^{-\frac{2 \alpha}{2+ l}},  \quad \lambda \to \infty,
\end{equation*}
for all $\varphi \in C_0(\mathbb R^N)$ satisfying $\varphi \ge 0$ and $\varphi(0) \ne 0$.}
\end{rem}
\begin{To}\label{thmupper2bis}
Let $N \geq 1$, $\alpha > 0$, $0 < l < N\alpha$, and let $C_l$ be as defined in \eqref{constanttype}. Let $\frac{l}{\alpha}\leq \gamma< \min \{N, \frac{2+ l}{\alpha}\}$ and let $\omega \in L^\infty(\mathbb R^N)$ be homogeneous of degree $0$, $\omega \geq 0$  and $\omega \not \equiv 0.$ Assume that $\varphi \in L^1_{\mathrm{loc}}(\mathbb{R}^N)$, $\varphi \geq 0$, satisfies
\begin{equation*}
 \varphi(x) \geq \omega(x)\,|x|^{-\gamma}, \quad \text{for } |x| \leq \varepsilon, 
\end{equation*}
for some $\varepsilon > 0$. Then, $T_{\max}(\lambda \varphi) < \infty$ for all sufficiently large $\lambda > 0$, and 
\begin{equation*}
  \limsup_{\lambda \to \infty} \lambda^{\left(\frac{2+ l}{2 \alpha}-\frac{\gamma}{2} \right)^{-1}} T_{\max}(\lambda \varphi) \le  \left(\frac{2+l}{2 \alpha}\right)^{\left(\frac{2+ l}{2}-\frac{\alpha \gamma}{2}\right)^{-1}} \left[\frac{C_l}{ \left\|e^{\Delta}(\omega |\cdot|^{-\gamma}) \right\|_{L^\infty}}\right]^{\left(\frac{2+ l}{2 \alpha}- \frac{\gamma}{2}\right)^{-1}}.
\end{equation*}
\end{To}
\begin{rem}$\;${\rm
\begin{itemize}
  \item [1)] Let $N \geq 1$, $\alpha > 0$ and  $0 < l < N\alpha$. Let $\frac{l}{\alpha}<\gamma< \min \{N, \frac{2+ l}{\alpha}\}$ and let $\varphi$ as in Theorem \ref{thmupper2bis} such that  $\varphi \in L_{\gamma}^{\infty}(\mathbb{R}^N)$. Then,  Theorem \ref{thm lower bound 2} and Theorem \ref{thmupper2bis} together show that
\begin{equation*}
  T_{\max}(\lambda \varphi) \sim \lambda^{-(\frac{2+ l}{2 \alpha}-\frac{\gamma}{2})^{-1}}, \quad  \lambda \to \infty.
\end{equation*}
  \item [2)] Let $N \geq 1$, $\alpha > 0$ and  $0 < l < N\alpha$. Let $\gamma= \frac{l}{\alpha}$   and let $\varphi$ as in Theorem \ref{thmupper2bis} such that  $\varphi \in L_{\frac{l}{\alpha}}^{\infty}(\mathbb{R}^N)$. Then,  Theorem \ref{thm lower bound 2bis} and Theorem \ref{thmupper2bis} together show that
\begin{equation*}
  T_{\max}(\lambda \varphi) \sim \lambda^{- \alpha}, \quad  \lambda \to \infty,
\end{equation*}
in agreement with the equivalence \eqref{Llinfty}.
  \item [3)] Proposition C.3 in \cite{TW4} provides an upper bound with the same rate in $\lambda$, but this result holds only in the subcritical case, that is, when $\alpha < \frac{2 + l}{N}$, and further assumes that $\gamma> \frac{l}{\alpha}$.
\end{itemize}}
\end{rem}
\begin{To}\label{thmupper2}
Let $N \geq 1$, $\alpha > 0$,  $-\min\{2, N\}<l<0$, and let $C_l$ be given by \eqref{constanttype}. Let $0<\gamma< \min \{N, \frac{2+ l}{\alpha}\}$ and let $\omega \in L^\infty(\mathbb R^N)$ be homogeneous of degree $0$, $\omega \geq 0$  and $\omega \not \equiv 0.$ Assume that $\varphi \in L^1_{\mathrm{loc}}(\mathbb{R}^N)$, $\varphi \geq 0$, satisfies
\begin{equation*}
 \varphi(x) \geq \omega(x)\,|x|^{-\gamma}, \quad \text{for } |x| \leq \varepsilon, 
\end{equation*}
for some $\varepsilon > 0$. Then, $T_{\max}(\lambda \varphi) < \infty$ for all sufficiently large $\lambda > 0$, and 
\begin{equation*}
  \limsup_{\lambda \to \infty} \lambda^{\left(\frac{2+ l}{2 \alpha}-\frac{\gamma}{2} \right)^{-1}} T_{\max}(\lambda \varphi) \le  \left(\frac{2+l}{2 \alpha}\right)^{\left(\frac{2+ l}{2}- \frac{\alpha \gamma}{2}\right)^{-1}} \left[\frac{C_l}{e^{\Delta}(\omega |\cdot|^{-\gamma})(0)}\right]^{\left(\frac{2+ l}{2 \alpha}- \frac{\gamma}{2}\right)^{-1}}.
\end{equation*}
\end{To}
\begin{rem}$\;${\rm
\begin{itemize}
  \item [1)] Let $N \geq 1$, $\alpha > 0$ and $-\min\{2, N\}<l<0$. Let $0<\gamma< \min \{N, \frac{2+ l}{\alpha}\}$ and let $\varphi$ as in Theorem \ref{thmupper2} such that  $\varphi \in L_{\gamma}^{\infty}(\mathbb{R}^N)$. Then, it follows from Theorems \ref{thm lower bound 2} and \ref{thmupper2} that
\begin{equation*}
  T_{\max}(\lambda \varphi) \sim \lambda^{-(\frac{2+ l}{2 \alpha}-\frac{\gamma}{2})^{-1}}, \quad  \lambda \to \infty.
\end{equation*}
  \item [2)] An upper bound with the same rate in $\lambda$ is obtained in Proposition C.3 of \cite{TW4}; however, it holds only in the subcritical case.
\end{itemize}}
\end{rem}
We now turn to the case of small $\lambda$. The following results hold.
\begin{To}\label{thmupper3bis}
Let $N \geq 1$, $\alpha > 0$, $0 < l < N\alpha$, and let $C_l$ be as defined in \eqref{constanttype}. Let $\frac{l}{\alpha} \le \gamma< \min \{N, \frac{2+ l}{\alpha}\}$ and let $\omega \in L^\infty(\mathbb R^N)$ be homogeneous of degree $0$, $\omega \geq 0$  and $\omega \not \equiv 0.$ Assume that $\varphi \in L^1_{\mathrm{loc}}(\mathbb{R}^N)$, $\varphi \geq 0$, satisfies
\begin{equation*}
 \varphi(x) \geq \omega(x)\,|x|^{-\gamma}, \quad \text{for } |x| \geq R, 
\end{equation*}
for some $R > 0$. Then, $T_{\max}(\lambda \varphi) < \infty$ for all  $\lambda > 0$, and 
\begin{equation*}
  \limsup_{\lambda \to 0} \lambda^{\left(\frac{2+ l}{2 \alpha}-\frac{\gamma}{2} \right)^{-1}} T_{\max}(\lambda \varphi) \le  \left(\frac{2+l}{2 \alpha}\right)^{\left(\frac{2+ l}{2}-\frac{\alpha \gamma}{2}\right)^{-1}} \left[\frac{C_l}{\left\|e^{\Delta}(\omega |\cdot|^{-\gamma})\right\|_{L^\infty}}\right]^{\left(\frac{2+ l}{2 \alpha}- \frac{\gamma}{2}\right)^{-1}}.
\end{equation*}
\end{To}
\begin{rem}$\;${\rm
\begin{itemize}
  \item [1)] Let $N \geq 1$, $\alpha > 0$ and  $0 < l < N\alpha$. Let $\frac{l}{\alpha}<\gamma< \min \{N, \frac{2+ l}{\alpha}\}$ and let $\varphi$ as in Theorem \ref{thmupper3bis} such that  $\varphi \in L_{\gamma}^{\infty}(\mathbb{R}^N)$. Then, by combining Theorems \ref{thm lower bound 2} and \ref{thmupper3bis}, we deduce that
\begin{equation*}
  T_{\max}(\lambda \varphi) \sim \lambda^{-(\frac{2+ l}{2 \alpha}-\frac{\gamma}{2})^{-1}}, \quad  \lambda \to 0.
\end{equation*}
  \item [2)]  Let $N \geq 1$, $\alpha > 0$ and  $0 < l < N\alpha$. Set $\gamma= \frac{l}{\alpha}$,   and let $\varphi$ as in Theorem \ref{thmupper3bis}, with  $\varphi \in L_{\frac{l}{\alpha}}^{\infty}(\mathbb{R}^N)$. Then, Theorems \ref{thm lower bound 2bis} and \ref{thmupper3bis} together imply that
\begin{equation*}
  T_{\max}(\lambda \varphi) \sim \lambda^{- \alpha}, \quad  \lambda \to 0.
\end{equation*}
  \item [3)] Proposition C.3 in \cite{TW4} establishes an upper bound with the same rate in $\lambda$, but this result is valid only in the subcritical case, and under the additional condition $\gamma> \frac{l}{\alpha}$.
  \item [4)] In \cite[Remark 5.2]{CFL}, the authors established that, for $\varphi=(1+|\cdot|^2)^{-\frac{\gamma}{2}}$ and sufficiently small $\lambda>0$, the maximal existence time $T_{\max}(\lambda\varphi)$ satisfies $\lambda^{\left(\frac{2+l}{2\alpha}-\frac{\gamma}{2}\right)^{-1}}T_{\max}(\lambda\varphi)\le C$, where $C>0$ is a constant whose explicit value was not determined.
\end{itemize}}
\end{rem}
\begin{To}\label{thmupper3}
Let $N \geq 1$, $\alpha > 0$,  $-\min\{2, N\}<l<0$, and let $C_l$ be given by \eqref{constanttype}.  Let $0<\gamma< \min \{N, \frac{2+ l}{\alpha}\}$ and let $\omega \in L^\infty(\mathbb R^N)$ be homogeneous of degree $0$, $\omega \geq 0$  and $\omega \not \equiv 0.$ Assume that $\varphi \in L^1_{\mathrm{loc}}(\mathbb{R}^N)$, $\varphi \geq 0$, satisfies
\begin{equation*}
 \varphi(x) \geq \omega(x)\,|x|^{-\gamma}, \quad \text{for } |x| \geq R, 
\end{equation*}
for some $R > 0$. Then, $T_{\max}(\lambda \varphi) < \infty$ for all  $\lambda > 0$, and 
\begin{equation*}
  \limsup_{\lambda \to 0} \lambda^{\left(\frac{2+ l}{2 \alpha}-\frac{\gamma }{2}\right)^{-1}} T_{\max}(\lambda \varphi) \le  \left(\frac{2+l}{2 \alpha}\right)^{\left(\frac{2+ l}{2}-\frac{\alpha \gamma}{2}\right)^{-1}} \left[\frac{C_l}{e^{\Delta}(\omega |\cdot|^{-\gamma})(0)}\right]^{\left(\frac{2+ l}{2 \alpha}- \frac{\gamma}{2}\right)^{-1}}.
\end{equation*}
\end{To}
\begin{rem}$\;${\rm
\begin{itemize}
  \item [1)] Let $N \geq 1$, $\alpha > 0$ and $-\min\{2, N\}<l<0$. Let $0<\gamma< \min \{N, \frac{2+ l}{\alpha}\}$ and let $\varphi$ as in Theorem \ref{thmupper3} be such that  $\varphi \in L_{\gamma}^{\infty}(\mathbb{R}^N)$. Then, Theorems \ref{thm lower bound 2} and \ref{thmupper3} together imply that
\begin{equation*}
  T_{\max}(\lambda \varphi) \sim \lambda^{-(\frac{2+ l}{2 \alpha}-\frac{\gamma}{2})^{-1}}, \quad  \lambda \to 0.
\end{equation*}
  \item [2)] Proposition C.3 in \cite{TW4} provides an upper bound with the same rate for $\lambda$, but this result is valid only in the subcritical case.
  \item [3)] In \cite[Remark 5.2]{CFL}, the authors showed that $\lambda^{\left(\frac{2+l}{2\alpha}-\frac{\gamma}{2}\right)^{-1}}T_{\max}(\lambda\varphi)\le C$ for sufficiently small $\lambda>0$, and $\varphi=(1+|\cdot|^2)^{-\frac{\gamma}{2}}.$  $C>0$ is a constant not specified.
\end{itemize}}
\end{rem}

In proving our results, we draw on the ideas of \cite{CFL}, \cite{TW4}, and \cite{TW5}. However, these arguments cannot be directly adapted to our setting, and new techniques are needed to establish the blow-up results and derive the upper bound estimates.

The remainder of the paper is organized as follows. In Section 2, we establish the necessary conditions for the local existence. In Section 3, we prove the blow-up criteria stated in Theorems \ref{thmblowupHénon} and \ref{thmblowup}. Finally, Section 4 is devoted to establishing the upper bound estimates stated in Theorems \ref{thmupper1bis}--\ref{thmupper3}.

Throughout this paper, for positive functions $f$ and $g$, we write $f(x) \sim g(x)$ as $x \to x_0$ if there exist two positive constants $C_1$ and $C_2$ such that $C_1  g(x) \leq f(x) \leq C_2 g(x)$, for all $x$ in a neighborhood of $x_0$.

%%%%%%%%%%%%%%%%%%%%%%%%%%%%%%%%%%%%%%%%%%%%%%%%%%%%%%%%%%%%%%%%%%%%%%%%%
\section{Necessary conditions for the existence}
\setcounter{equation}{0}
%%%%%%%%%%%%%%%%%%%%%%%%%%%%%%%%%%%%%%%%
\label{nececarycond}
In this section, we give the necessary conditions for the existence of a nonnegative solution to \eqref{NLHHint}. Similar conditions were obtained in \cite[Theorem 1, p. 546]{W4} for the case $l = 0$, that is, for the classical nonlinear heat equation. We begin by recalling the following lemma from \cite[Theorem 2.12, p. 27]{AF}.
\begin{lem}[A Reverse H\"older inequality]\label{Reverse Holder}
Let $\Omega$ be a subset of $\mathbb{R}^N$ and let $p \in (1,\infty)$. Suppose that $f, g \in L^{1}(\Omega)$ and $|g|>0$ a.e. in $\Omega$. Then,
\begin{equation*}
\int_{\Omega} |f(y)|^p |g(y)|^{1-p} \, dy \geq  
\left( \int_{\Omega} |f(y)| \, dy \right)^p 
\left( \int_{\Omega} |g(y)| \, dy \right)^{1-p}.
\end{equation*}
\end{lem}
We also recall  the commutation relation between $D_\lambda$ and the heat semigroup on $\mathbb{R}^N$.
\begin{equation}\label{commutation}
D_\lambda e^{\tau\lambda^2\Delta} = e^{\tau\Delta}D_\lambda, \quad \forall  \lambda > 0, \quad \tau > 0.
\end{equation}
In particular, we have
\begin{equation}\label{eq45}
   D_{\sqrt t} e^{t \Delta} = e^{\Delta} D_{\sqrt t}, \quad \forall t>0.
\end{equation}
We have the following result.
\begin{prop}\label{localexistence}
Let $N \geq 1$, $\alpha > 0$, $0 \leq l < N\alpha$ and $C_l$ the constant defined in \eqref{constanttype}. Let $u_0 \geq 0$ be a locally integrable function  on $\mathbb R^N$, and let $T_{\max}(u_0)$ denote the maximal existence time of a nonnegative solution $u$ to \eqref{NLHHint} on $(0, T_{\max}(u_0)) \times \mathbb{R}^N$.  Then, 
\begin{itemize}
  \item [(a)] for all $(t,x) \in (0,T_{\max}(u_0)) \times \mathbb{R}^N$, we have
\begin{equation}\label{neccd}
\frac{2}{2+l} \alpha  t^{\frac{2+l}{2}} \left[\left(e^{t\Delta}u_0\right)(x)\right]^\alpha \leq C_l^\alpha.
\end{equation}
Consequently,
\begin{equation}\label{neccdinfty}
\frac{2}{2+l} \alpha  t^{\frac{2+l}{2}} \|e^{t\Delta}u_0\|_{L^\infty}^\alpha \leq C_l^\alpha,
\end{equation}
for all $t \in (0,T_{\max}(u_0))$.
  \item [(b)] If $0<T_{\max}(u_0)<\infty$, then
  \begin{equation}\label{neccdinftyinfty}
\frac{2}{2+l} \alpha  \left(T_{\max}(u_0)\right)^{\frac{2+l}{2}} \|e^{T_{\max}(u_0)\Delta}u_0\|_{L^\infty}^\alpha \leq C_l^\alpha.
\end{equation}
\item [(c)] If $0<T_{\max}(u_0)<\infty$, $u_0=\lambda \varphi$ with $\lambda>0$, then
  \begin{equation}\label{neccdinftyinftybis}
\lambda^\alpha \frac{2}{2+l} \alpha  \left(T_{\max}(\lambda \varphi)\right)^{\frac{2+l}{2}} \|e^{T_{\max}(\lambda \varphi)\Delta}\varphi\|_{L^\infty}^\alpha \leq C_l^\alpha.
\end{equation}
\end{itemize}
\end{prop}
\begin{proof}
\begin{itemize}
  \item [(a)] The case $l = 0$ has been treated in \cite[Theorem 1, p. 546]{W4}, whereas the case $0 < l < N\alpha$ is treated in \cite{CFL}.  See also \cite[Lemma 4.1]{Tsusaka}.
  \item [(b)] Since inequality \eqref{neccd} is valid for every $0 < t < T_{\max}(u_0)$, it follows that it remains valid at $t = T_{\max}(u_0)$.
  \item [(c)] Inequality \eqref{neccdinftyinftybis} is an immediate consequence of \eqref{neccdinftyinfty}.
\end{itemize}
\end{proof}
In the Hénon case we have the following result using weighted spaces.
\begin{prop}\label{localexistencebis}
Let $N \geq 1$, $\alpha > 0$ and  $0 \leq l < N\alpha$. Let $u_0 \geq 0$ be a locally integrable function  on $\mathbb R^N$, and let $T_{\max}(u_0)$ denote the maximal existence time of a nonnegative solution $u$ to \eqref{NLHHint} on $(0, T_{\max}(u_0)) \times \mathbb{R}^N$.  Then,
\begin{itemize}
  \item [(a)] for every $(t,x) \in (0,T_{\max}(u_0)) \times \mathbb{R}^N$, we have
\begin{equation}\label{neccdz}
\alpha  t \left[ |x|^{\frac{l}{\alpha}} \left(e^{t\Delta}u_0\right)(x)\right]^\alpha  \leq  \| {e}^{\Delta} |\cdot|^{-\frac{l}{\alpha}} \|_{L_{\frac{l}{\alpha}}^\infty}^{\alpha}.
\end{equation}
Therefore, for each $t \in (0,T_{\max}(u_0))$, 
\begin{equation}\label{neccdinftyz}
\alpha  t \|e^{t\Delta}u_0\|_{L_{\frac{l}{\alpha}}^\infty}^\alpha \leq   \| {e}^{\Delta} |\cdot|^{-\frac{l}{\alpha}} \|_{L_{\frac{l}{\alpha}}^\infty}^{\alpha}.
\end{equation}
  \item [(b)] If $0<T_{\max}(u_0)<\infty$, then
  \begin{equation}\label{neccdinftyinftyz}
   \alpha  T_{\max}(u_0) \|e^{T_{\max}(u_0)\Delta}u_0\|_{L_{\frac{l}{\alpha}}^\infty}^\alpha \leq  \| {e}^{\Delta} |\cdot|^{-\frac{l}{\alpha}} \|_{L_{\frac{l}{\alpha}}^\infty}^{\alpha}.
  \end{equation}
\item [(c)] If $0<T_{\max}(u_0)<\infty$, $u_0=\lambda \varphi$ with $\lambda>0$, then
  \begin{equation}\label{neccdinftyinftybisz}
\lambda^\alpha \alpha T_{\max}(\lambda \varphi) \|e^{T_{\max}(\lambda \varphi)\Delta}\varphi\|_{L_{\frac{l}{\alpha}}^\infty}^\alpha \leq \| {e}^{\Delta} |\cdot|^{-\frac{l}{\alpha}} \|_{L_{\frac{l}{\alpha}}^\infty}^{\alpha}.
\end{equation}
\end{itemize}
\end{prop}
\begin{proof}
\begin{itemize}
  \item [(a)]The proof relies on several arguments from \cite{CFL, Tayachi1, W4}. Without loss of generality, we assume that $u_0 \not\equiv 0$ and $0 < l < N\alpha$. Fix $\tau \in (0,T_{\max}(u_0))$, and let $0<t<\tau$. By Fubini's theorem, we have
\begin{equation}\label{eq124}
   e^{(\tau-t) \Delta} u(t)= e^{\tau \Delta}u_0+ \int_0^t {e}^{(\tau-\sigma)\Delta} \big[|\cdot|^l u^{\alpha+1}(\sigma)\big] d\sigma.
 \end{equation}
 Using Lemma \ref{Reverse Holder} (A reverse Hölder inequality) with $f(y)=G(\tau-\sigma, x-y) u(\sigma,y)$, $g(y)=G(\tau-\sigma, x-y)|y|^{-\frac{l}{\alpha}}$ and $p=\alpha+1$, we obtain
\begin{eqnarray}\label{eqn145}
    {e}^{(\tau-\sigma)\Delta} \big[|\cdot|^l u^{\alpha+1}(\sigma)\big](x) &=&  \int_{\mathbb R^N} G(\tau-\sigma, x-y) |y|^l u^{\alpha+1}(\sigma, y) dy \nonumber\\
    &=&  \int_{\mathbb R^N} \left( G(\tau-\sigma, x-y) |y|^{-\frac{l}{\alpha}} \right)^{-\alpha} \times  \nonumber\\&&  \hspace{2 cm}\left( G(\tau-\sigma, x-y) u(\sigma, y) \right)^{\alpha+1} dy \nonumber\\
    &\geq& \left(\int_{\mathbb R^N}  G(\tau-\sigma, x-y) |y|^{-\frac{l}{\alpha}}   dy\right)^{-\alpha} \times \nonumber\\ && \hspace{2 cm}\left(\int_{\mathbb R^N}   G(\tau-\sigma, x-y) u(\sigma, y)  dy \right)^{\alpha+1} \nonumber\\
    &=&   \left({e}^{(\tau-\sigma)\Delta} |\cdot|^{-\frac{l}{\alpha}}\right)^{-\alpha}(x)  \left({e}^{(\tau-\sigma)\Delta} u(\sigma) \right)^{\alpha+1}(x),
 \end{eqnarray}
 for all $x\in \mathbb R^N$ and $0<\sigma<\tau$. Combining \eqref{eq124} and \eqref{eqn145}, we get
\begin{equation*}
    e^{(\tau-t) \Delta} u(t) \geq  e^{\tau \Delta}u_0+ \int_0^t \left({e}^{(\tau-\sigma)\Delta} |\cdot|^{-\frac{l}{\alpha}}\right)^{-\alpha}  \left({e}^{(\tau-\sigma)\Delta} u(\sigma) \right)^{\alpha+1} d\sigma:= \Phi(t).
 \end{equation*}
 We have,
 \begin{eqnarray*}
   \Phi'(t) &=& \left({e}^{(\tau- t)\Delta} |\cdot|^{-\frac{l}{\alpha}}\right)^{-\alpha}  \left({e}^{(\tau-t)\Delta} u(t) \right)^{\alpha+1} \\
    &\geq&  \left({e}^{(\tau- t)\Delta} |\cdot|^{-\frac{l}{\alpha}}\right)^{-\alpha} \left(\Phi(t) \right)^{\alpha+1}.
 \end{eqnarray*}
Therefore,
\begin{equation*}
   \left(\Phi^{-\alpha}(t)\right)'=-\alpha  \Phi'(t) \left(\Phi(t) \right)^{-\alpha-1} \leq -\alpha \left({e}^{(\tau- t)\Delta} |\cdot|^{-\frac{l}{\alpha}}\right)^{-\alpha}.
 \end{equation*}
By integrating both sides of the above inequality over $(0, \tau)$, we obtain
\begin{equation}\label{ineq1452}
  \Phi^{-\alpha}(\tau)- \Phi^{-\alpha}(0) \leq -\alpha  \int_0^{\tau}\left({e}^{(\tau- t)\Delta} |\cdot|^{-\frac{l}{\alpha}}\right)^{-\alpha} dt.
\end{equation}
Since $\Phi(0) = e^{\tau \Delta} u_0$ and $\Phi(\tau) > 0$, it follows from inequality \eqref{ineq1452} that
\begin{equation}\label{eq 47}
 \left(e^{\tau \Delta}u_0\right)^\alpha  \int_0^{\tau}\left({e}^{(\tau- t)\Delta} |\cdot|^{-\frac{l}{\alpha}}\right)^{-\alpha} dt \leq \frac{1}{\alpha}.
\end{equation}
Inequality \eqref{eq 47} can be rewritten as
\begin{equation}\label{eq 47bis}
 \left( |\cdot|^{\frac{l}{\alpha}} e^{\tau \Delta}u_0\right)^\alpha  \int_0^{\tau}\left( |\cdot|^{\frac{l}{\alpha}} {e}^{(\tau- t)\Delta} |\cdot|^{-\frac{l}{\alpha}}\right)^{-\alpha} dt \leq \frac{1}{\alpha}.
\end{equation}
Consequently, \eqref{eq 47bis} yields
\begin{equation}\label{eq 874}
    \left( |\cdot|^{\frac{l}{\alpha}} e^{\tau \Delta}u_0\right)^\alpha \int_0^{\tau} \left \| {e}^{(\tau- t)\Delta} |\cdot|^{-\frac{l}{\alpha}} \right \|_{L_{\frac{l}{\alpha}}^\infty}^{-\alpha} dt \leq \frac{1}{\alpha}.
\end{equation}
From \eqref{eq45} and the homogeneity of $|\cdot|^{-\frac{l}{\alpha}}$, we have
\begin{eqnarray}\label{eq 417l}
    \left\|e^{\Delta}|\cdot|^{- \frac{l}{\alpha}}\right\|_{L_{\frac{l}{\alpha}}^\infty}=  \left\|D_{\sqrt{t}} e^{ t \Delta} D_{\frac{1}{\sqrt{t}}}|\cdot|^{-\frac{l}{\alpha}}\right\|_{L_{\frac{l}{\alpha}}^\infty}&=& t^{-\frac{l}{2 \alpha}} \left\|e^{ t \Delta} D_{\frac{1}{\sqrt{t}}}|\cdot|^{- \frac{l}{\alpha}}\right\|_{L_{\frac{l}{\alpha}}^\infty} \nonumber\\
   &=& \left\|e^{ t \Delta} |\cdot|^{-\frac{l}{\alpha}}\right\|_{L_{\frac{l}{\alpha}}^\infty}, \quad \forall t>0.
\end{eqnarray}
Therefore, \eqref{neccdz} follows from \eqref{eq 874} and \eqref{eq 417l}.
  \item [(b)] Since inequality \eqref{neccdz} holds for every $0 < t < T_{\max}(u_0)$, it also holds at $t = T_{\max}(u_0)$.
  \item [(c)] \eqref{neccdinftyinftybisz} follows immediately from \eqref{neccdinftyinftyz}.
\end{itemize}
\end{proof}
\begin{rem}\label{rem 124}
{\rm For $0 < l < N\alpha$, Lemma 2.1 of \cite{CIT} implies that   $\|e^{\Delta} |\cdot|^{-\frac{l}{\alpha}} \|_{L_{\frac{l}{\alpha}}^\infty} \le C \left\||\cdot|^{-\frac{l}{\alpha}} \right\|_{L_{\frac{l}{\alpha}}^\infty}=C$.}
\end{rem}
 \begin{prop}\label{localexistencelneg}
Let $N\geq 1,$ $\alpha>0$, $-\min\{2, N\}<l<0$ and $C_l$ the constant defined in \eqref{constanttype}. Let $ u_0 \geq 0$  be a locally integrable function  on $\mathbb R^N$, and let $T_{\max}(u_0)$ denote the maximal existence time of a nonnegative solution $u$ of  \eqref{NLHHint} on $\left(0, T_{\max}(u_0) \right)\times \mathbb R^N$.  
Then,
\begin{itemize}
  \item [(a)] for all $t \in (0, T_{\max}(u_0))$, we have
\begin{equation}
\label{neccdbis}
{2\over 2+l} \alpha t^{{2+l\over 2}} \left[({e}^{t\Delta}u_0)(0)\right]^\alpha \le C_l^\alpha.
\end{equation}
  \item [(b)]  If $0<T_{\max}(u_0)<\infty$, then
  \begin{equation}\label{neccdbisbis}
{2\over 2+l} \alpha \left(T_{\max}(u_0) \right)^{{2+l\over 2}} \left[({e}^{ T_{\max}(u_0) \Delta}u_0)(0)\right]^\alpha \le C_l^\alpha.
\end{equation}
  \item [(c)]  If $0<T_{\max}(u_0)<\infty$, $u_0=\lambda \varphi$ with $\lambda>0$, then
  \begin{equation}\label{neccdbisbisbis}
\lambda^{\alpha}{2\over 2+l} \alpha \left(T_{\max}(\lambda \varphi) \right)^{{2+l\over 2}} \left[({e}^{ T_{\max}(\lambda \varphi) \Delta} \varphi)(0)\right]^\alpha \le C_l^\alpha.
\end{equation}
\end{itemize}
\end{prop}
\begin{proof}
\begin{itemize}
  \item [(a)] The proof is given in \cite{CFL} and relies on inequality \eqref{eq 47} together with the identity  $\left({e}^{(\tau- t)\Delta} |\cdot|^{-\frac{l}{\alpha}}\right)(0)= C_l (\tau-t)^{-\frac{l}{2 \alpha}},$ where $C_l$ is defined by \eqref{constanttype}. See also \cite[Lemma 4.1]{Tsusaka}.
  \item [(b)] Inequality \eqref{neccdbis} is valid for all $0 < t < T_{\max}(u_0)$, and hence remains valid at $t = T_{\max}(u_0)$.
  \item [(c)] \eqref{neccdbisbisbis} follows directly from \eqref{neccdbisbis}.
\end{itemize}
\end{proof}
 
\begin{rem}
{\rm For $-\min\{2, N\} < l < 0$. We cannot use the bound $e^{t \Delta}(|\cdot|^{-\frac{l}{\alpha}}) \leq e^{t \Delta}(|\cdot|^{-\frac{l}{\alpha}})(0)$ to obtain more general estimate analogous to \eqref{neccd}. For example, when $N=1$, $\alpha= \frac{1}{4}$ and $l=- \frac{1}{2}$, it is not difficult to see that $e^{t \Delta}(|\cdot|^{-\frac{l}{\alpha}})(x)= x^2+ 2 t \geq e^{t \Delta}(|\cdot|^{-\frac{l}{\alpha}})(0)$, for all $(t, x)\in (0, \infty) \times \mathbb R$.} 
\end{rem}
\section{Blow-up results}
Let $0<\gamma<N$. We define the Banach space ${\mathcal X}_\gamma$,  introduced in \cite{CDW-DCDS}, see also  \cite{MTW, TW5},  by
\begin{equation*}
  {\mathcal X}_\gamma:= \left\{ \psi \in L^1_{loc}(\mathbb R^N);\quad   {|\cdot|^{\gamma}\psi} \in L^{\infty}(\mathbb R^N) \right\}=L^\infty_\gamma,
\end{equation*}
with norm $\|\psi\|_{\mathcal{X}_\gamma} := \left\| |\cdot|^{\gamma}\psi \right\|_{L^\infty},$ for all $\psi \in \mathcal{X}_\gamma.$ Hence, $|\psi(x)| \leq  \|\psi\|_{{\mathcal X}_\gamma}  |x|^{-\gamma}$ for almost every  $x \in \mathbb R^N$. We also consider the Banach space
\begin{equation}\label{ygamma}
  \mathcal Y_\gamma:= \left\{ \psi \in L^1_{loc}(\mathbb R^N);\quad |\cdot|^{-\gamma} \psi \in L^{1}(\mathbb R^N) \right\}.
\end{equation}
It is straightforward to see that the predual of $\mathcal{X}\gamma$ is  $\mathcal{Y}\gamma$. Let $K>0$, and define the closed ball  
\begin{equation}\label{BK}
  \mathcal B_K:= \left\{\psi \in \mathcal X_\gamma; \quad \|\psi\|_{\mathcal X_\gamma}\leq K \right\}.
\end{equation}
Let $\mathcal{B}_K^*$ denote $\mathcal{B}_K$ endowed with the weak-$^*$ topology on $\X_\gamma.$  It is well known that $\mathcal{B}_K^*$ is a compact metric space (see \cite{brezis}), and therefore complete and separable. Moreover, by \cite[Proposition 2.1, p. 1110]{CDW-DCDS}, if $(\psi_n)_{n=1}^\infty$ is a sequence in $\mathcal{B}_K$ and $\psi \in \mathcal{B}_K$, then $\psi_n \to \psi$ in $\mathcal{B}_K^*$ if and only if $\psi_n \to \psi$ in $\mathcal{D}'(\mathbb{R}^N)$. For $\psi\in {\mathcal B}_K$, we recall from \cite[Formula (1.17), p. 1107]{CDW-DCDS} the  following set
\begin{equation}
\label{Zu0}
\mathcal{Z}^\gamma(\psi) := \left\{ z \in \mathcal{B}_K;\quad \exists\, \lambda_n \to \infty \text{ such that } \lambda_n^\gamma D_{\lambda_n} \psi \to z \text{ in } \mathcal{B}_K^* \right\}, 
\end{equation}
where $D_\lambda$ is the dilation operator defined in \eqref{dilation}. By \cite[Remark 2.4 (iii), p. 1112]{CDW-DCDS}, $\mathcal{Z}^\gamma(\psi)$ is a nonempty, compact, and connected subset of the compact metric space $\mathcal{B}_K^*$. 

The proofs of Theorems \ref{thmblowupHénon} and \ref{thmblowup}, concerning blow-up for equation \eqref{NLHH}, are entirely based on properties of the linear heat semigroup. The main ingredients are  Propositions \ref{localexistence} and \ref{localexistencelneg}, which provide necessary conditions on the heat semigroup for the existence of nonnegative solutions on a given time interval, and  \cite[Proposition 3.8, p. 1123]{CDW-DCDS}. The following lemma will be needed. Although a proof can be found in \cite[Exercise 4.26, p. 126]{brezis}, we present a different one here.
\begin{lem}\label{lem3}
Let $(u_n)_{n \geq 0}$ be a sequence of functions in $L^1_{\mathrm{loc}}(\mathbb{R}^N)$. Suppose that
\begin{itemize}
  \item[(i)] $u_n \geq 0$ a.e. for all $n \geq 0$,
  \item[(ii)] there exists $u \in L^1_{\mathrm{loc}}(\mathbb{R}^N)$ such that $u_n \to u$ in $\mathcal{D}'(\mathbb{R}^N)$.
\end{itemize}
Then $u \geq 0$ a.e. in $\mathbb{R}^N$.
\end{lem}
\begin{proof}
Since $u_n \to u$ in $\mathcal{D}'(\mathbb{R}^N)$, we have
\begin{equation*}
\int_{\mathbb{R}^N} u_n(x) \varphi(x) dx \longrightarrow \int_{\mathbb{R}^N} u(x) \varphi(x) dx, \quad \text{for all } \varphi \in C_c^\infty(\mathbb{R}^N).
\end{equation*}
Using the fact that $u_n \geq 0$ a.e., it follows that for all  $\varphi \in C_c^\infty(\mathbb{R}^N),$ $\varphi \geq 0$,
\begin{equation}\label{1278h}
\int_{\mathbb{R}^N} u(x) \varphi(x) dx \geq 0.
\end{equation}
In particular, \eqref{1278h} implies that 
\begin{equation}\label{compact1}
\int_K u(x) dx \geq 0, \quad \text{for every compact set } K \subset \mathbb{R}^N.
\end{equation}
See the proof of \cite[Corollary 4.24, p. 110]{brezis}. Now, suppose by contradiction that $A := \{ x \in \mathbb{R}^N;\quad u(x) < 0 \}$ has positive measure. By the inner regularity of the Lebesgue measure \cite[Proposition 7.5, p. 216]{Folland}, there exists a compact set $K_1 \subset A$ with positive measure. Hence, 
\begin{equation*}
\int_{K_1} u(x) dx < 0,
\end{equation*}
which contradicts \eqref{compact1}.
\end{proof}
\begin{proof}[Proof of  Theorem \ref{thmblowupHénon}]
Let $u_0$ and $f$ be as in Theorem \ref{thmblowupHénon}. By comparison,  it is enough to show that $T_{\max}(f)<\infty$. We proceed by contradiction and assume that $T_{\max}(f) = \infty$. Using \eqref{eq45}, we obtain
\begin{equation*}
 \left\| e^{t \Delta}f \right\|_{L^\infty} =  \left\| D_{\sqrt t}e^{t \Delta}f\right\|_{L^\infty}= \left\| e^{ \Delta}  D_{\sqrt t}f \right\|_{L^\infty}, \quad \forall t>0.
\end{equation*}
Thus,
\begin{equation*}
    t^{\frac{\gamma}{2}}  \left \| e^{t \Delta}f \right\|_{L^\infty}  =\left\| e^{ \Delta}  \left( t^{\frac{\gamma}{2}} D_{\sqrt t}f \right) \right\|_{L^\infty},  \quad \forall t>0.
\end{equation*}
Let $z \in L^1_{loc}(\mathbb{R}^N)$ with $z \not\equiv 0$, and let $\lambda_n \to \infty$ be such that
\begin{equation*}
  \lambda_n^\gamma D_{\lambda_n} f \to z \quad \text{in } \mathcal{D}'(\mathbb{R}^N) \text{ as } n \to \infty.
\end{equation*}
By \cite[p. 169]{TW5}, it follows that $z \in \mathcal{Z}^\gamma(f)$ and that the convergence also holds in $\mathcal{B}_K^*$, where $K$ denotes the positive constant defined in \eqref{initialinfmodulef}.  Next, define $t_n := \lambda_n^2 \to \infty$ as $n \to \infty$. Then,
\begin{equation*}
  t_n^{\frac{\gamma}{2}} D_{\sqrt{t_n}} f \to z \quad \text{in } \mathcal{B}_K^*.
\end{equation*}
Applying \cite[Proposition 3.8 (i), p. 1123]{CDW-DCDS}, we get
\begin{equation*}
  e^\Delta \left( t_n^{\frac{\gamma}{2}} D_{\sqrt{t_n}} f \right) \to e^\Delta z \quad \text{in } C_0(\mathbb{R}^N).
\end{equation*}
Hence,
\begin{equation}\label{eq 74f}
\lim_{n \to \infty} t_n^{\frac{\gamma}{2}} \left\| e^{t_n \Delta}  f \right\|_{L^\infty} = \left\| e^\Delta z \right\|_{L^\infty}>0.
\end{equation}
From \eqref{neccdinfty}, we obtain
\begin{equation}\label{ineq 45h}
    t_n^{\frac{\gamma}{2}}  \left\| e^{t_n \Delta}  f \right\|_{L^\infty} \leq  C_l   \left(\frac{2+ l}{2 \alpha} \right)^{\frac{1}{\alpha}} t_n^{\frac{\gamma}{2}- \frac{2+ l}{2 \alpha}},  \quad \forall n>0.
\end{equation}
Passing to the limit as $n \to \infty$ in \eqref{ineq 45h} and using \eqref{eq 74f}, we obtain a contradiction when $\alpha < \frac{2+l}{\gamma}$, since $\left\| e^\Delta z \right\|_{L^\infty}> 0$. We also get a contradiction if $\alpha= \frac{2+ l}{\gamma}$ by letting $n \to \infty$, since in this case $\left\| e^\Delta z \right\|_{L^\infty} > C_l \left( \frac{2+l}{2 \alpha}\right)^{\frac{1}{\alpha}}$. Therefore, we conclude that $T_{\max}(f) < \infty$, which completes the proof. 
\end{proof}
\begin{proof}[Proof of  Theorem \ref{thmblowup}]
Let $u_0$ and $f$ be as in Theorem \ref{thmblowup}. It suffices to prove that $T_{\max}(f) < \infty$. We argue by contradiction and assume that $T_{\max}(f) = \infty$. It then follows from \eqref{eq45} that
\begin{equation*}
   t^{\frac{\gamma}{2}}  \left(e^{t \Delta}f \right)(0) =t^{\frac{\gamma}{2}}  \left(D_{\sqrt t}e^{t \Delta}f\right)(0)= t^{\frac{\gamma}{2}} \left(e^{ \Delta}  D_{\sqrt t}f\right)(0), \quad \forall t>0.
\end{equation*}
Thus, \eqref{neccdbis} implies that
\begin{equation}\label{ineq 1254}
   t^{\frac{\gamma}{2}}  \left(D_{\sqrt t}e^{t \Delta}f\right)(0)=  t^{\frac{\gamma}{2}}  \left(e^{t \Delta}f \right)(0)\leq  C_l   \left(\frac{2+ l}{2 \alpha} \right)^{\frac{1}{\alpha}} t^{\frac{\gamma}{2}- \frac{2+ l}{2 \alpha}},  \quad \forall t>0.
\end{equation}
Let $z \in L^1_{loc}(\mathbb{R}^N)$ with $z \not\equiv 0$, and let $\lambda_n \to \infty$ be such that
\begin{equation*}
  \lambda_n^\gamma D_{\lambda_n} f \to z \quad \text{in } \mathcal{D}'(\mathbb{R}^N) \text{ as } n \to \infty.
\end{equation*}
Let $t_n := \lambda_n^2$, so that $t_n \to \infty$ as $n \to \infty$. As in the proof of Theorem \ref{thmblowupHénon}, we obtain
\begin{equation*}
  e^\Delta \left( t_n^{\frac{\gamma}{2}} D_{\sqrt{t_n}} f \right) \to e^\Delta z \quad \text{in } C_0(\mathbb{R}^N),
\end{equation*}
and in particular,
\begin{equation*}
\lim_{n \to \infty} t_n^{\frac{\gamma}{2}} \left( e^\Delta D_{\sqrt{t_n}} f \right)(0) = \left( e^\Delta z \right)(0).
\end{equation*}
Using \eqref{subseqlimz}, Lemma \ref{lem3} and the fact that  $z \not\equiv 0$ a.e., we see that 
\begin{equation}\label{eq4784_restructured}
   \left( e^\Delta z \right)(0)>0.
\end{equation}
 Finally, replacing $t$ by $t_n$ in \eqref{ineq 1254} and taking the limit $n \to \infty$, we obtain, if $\alpha<{2+l\over \gamma}$,
\begin{equation*}
  (e^\Delta z)(0) \le 0,
\end{equation*}
which contradicts \eqref{eq4784_restructured}. We also get a contradiction if $\alpha= \frac{2+ l}{\gamma}$ by letting $n \to \infty$, since in this case $(e^{\Delta}z)(0) > C_l \left( \frac{2+l}{2 \alpha}\right)^{\frac{1}{\alpha}}$. Therefore, we conclude that $T_{\max}(f) < \infty$, which completes the proof. 
\end{proof}
\section{Upper bounds of the life-span}
In this section, we exploit conditions \eqref{neccd}, \eqref{neccdz} and \eqref{neccdbis} for the existence of solutions. Conditions  \eqref{neccd} and \eqref{neccdbis} imply that if $u$ is a nonnegative solution of \eqref{NLHHint} defined on its maximal existence interval $\left(0, T_{\max}(u_0)\right)$, then
\begin{equation}\label{cond 12}
\frac{2}{2+l}\,\alpha\, t^{\frac{2+l}{2}} \left[({e}^{t\Delta}u_0)(0)\right]^\alpha \le C_l^\alpha,
\end{equation}
for all $t \in \left(0, T_{\max}(u_0)\right)$, where $-\min\{2,N\} < l < N\alpha$ and $u_0 \geq 0$ is a locally integrable function on $\mathbb{R}^N$. In the particular case $u_0 = \lambda \varphi$, condition \eqref{cond 12} becomes
\begin{equation}\label{cond 12bis}
\lambda^\alpha \frac{2}{2+l} \alpha  t^{\frac{2+l}{2}} \left[({e}^{t\Delta}\varphi)(0)\right]^\alpha \le C_l^\alpha,
\end{equation}
for all $t \in (0, T_{\max}(\lambda \varphi))$. This estimate yields a blow-up criterion for large initial data. Fix $\varphi \geq 0$ with $\varphi \not\equiv 0$, and let $t>0$. If $\lambda > C_l \left(\frac{2+l}{2\alpha}\right)^{\frac{1}{\alpha}}   t^{-\frac{2+l}{2\alpha}} \left[({e}^{t\Delta}\varphi)(0)\right]^{-1},$ then \eqref{cond 12bis} is violated at time $t$. Consequently, $t$ cannot belong to the interval of existence, and therefore
\begin{equation*}
T_{\max}(\lambda \varphi) \le t.
\end{equation*}
Since $t>0$ is arbitrary, we deduce that for every $t>0$ there exists $A(t)>0$ such that $\lambda > A(t)$ implies $T_{\max}(\lambda \varphi) \le t$. This proves that for any nontrivial nonnegative function $\varphi$, one obtains
\begin{equation}\label{limtmax}
\lim_{\lambda \to \infty} T_{\max}(\lambda \varphi) = 0.
\end{equation}

We now give the proofs of the upper bounds.
\begin{proof}[Proof of Theorem \ref{thmupper1bis}]
Let $\varphi\in L_{\frac{l}{\alpha}}^\infty(\mathbb{R}^N)$ be such that $\varphi\geq 0$ and $\varphi\not\equiv 0$. By inequality \eqref{neccdinftyinftybisz}, we have $$\lambda^\alpha T_{\max}(\lambda \varphi) \le \frac{1}{\alpha} \left(\frac{ \| {e}^{\Delta} |\cdot|^{-\frac{l}{\alpha}} \|_{L_{\frac{l}{\alpha}}^\infty}}{\|e^{T_{\max}(\lambda \varphi)\Delta}\varphi\|_{L_{\frac{l}{\alpha}}^\infty}}\right)^{\alpha}.$$ Consequently, 
\begin{equation*}
  \limsup_{\lambda \to \infty} \lambda^\alpha T_{\max}(\lambda \varphi) \le \frac{1}{\alpha} \left(\frac{ \| {e}^{\Delta} |\cdot|^{-\frac{l}{\alpha}} \|_{L_{\frac{l}{\alpha}}^\infty}}{\liminf_{\lambda \to \infty} \|e^{T_{\max}(\lambda \varphi)\Delta}\varphi\|_{L_{\frac{l}{\alpha}}^\infty}}\right)^{\alpha}.
\end{equation*}
Since $T_{\max}(\lambda\varphi)\to0$ as $\lambda\to\infty$, it suffices to show that $\|\varphi\|_{L_{\frac{l}{\alpha}}^\infty} \le \liminf_{t \to 0} \|e^{t \Delta}\varphi\|_{L_{\frac{l}{\alpha}}^\infty}$. 

We claim that, for every $\psi \in  C_c^\infty(\mathbb R^N)$, 
\begin{itemize}
  \item [(i)] $\int_{\mathbb R^N} (e^{t \Delta}\varphi) \psi= \int_{\mathbb R^N} \varphi (e^{t \Delta}\psi)$. Indeed, by \eqref{semigroup1}, 
       \begin{equation*}
         \int_{\mathbb R^N} (e^{t \Delta}\varphi)(x) \psi(x) dx = \int_{\mathbb R^N} \left (\int_{\mathbb R^N} G(t, x-y) \varphi(y) dy \right) \psi(x) dx.
       \end{equation*}
       We first verify that the double integral is absolutely convergent.  In view of Remark \ref{rem 124}, and using the assumptions $0<l<N\alpha$ and $\psi\in C_c^\infty(\mathbb{R}^N)$, we obtain
      \begin{eqnarray*}
     \int_{\mathbb R^N} \int_{\mathbb R^N} G(t, x-y)  |\varphi(y)| |\psi(x)|  dy dx & \le&   \int_{\mathbb R^N} \int_{\mathbb R^N} G(t, x-y)  \|\varphi\|_{L_{\frac{l}{\alpha}}^\infty} |y|^{-\frac{l}{\alpha}}  |\psi(x)| dy dx  \\
     &=& \|\varphi\|_{L_{\frac{l}{\alpha}}^\infty}  \int_{\mathbb R^N}  \left(e^{t\Delta} |\cdot|^{-\frac{l}{\alpha}}\right)(x)  |\psi(x)|  dx \\
     &\le& \|\varphi\|_{L_{\frac{l}{\alpha}}^\infty}  \int_{\mathbb R^N} \|e^{t \Delta}  |\cdot|^{-\frac{l}{\alpha}} \|_{L_{\frac{l}{\alpha}}^\infty}   |x|^{-\frac{l}{\alpha}} |\psi(x)|  dx \\
     &\le&C  \|\varphi\|_{L_{\frac{l}{\alpha}}^\infty}  \int_{\mathbb R^N}    |x|^{-\frac{l}{\alpha}} |\psi(x)|  dx< \infty.  
  \end{eqnarray*}
  Hence, Fubini's theorem applies, and we can interchange the order of integration. Using the symmetry of the heat kernel, $G(t, x-y)=G(t, y-x)$, we obtain
  \begin{eqnarray*}
     \int_{\mathbb R^N} \left (\int_{\mathbb R^N} G(t, x-y) \varphi(y) dy \right) \psi(x) dx &=& \int_{\mathbb R^N} \left (\int_{\mathbb R^N} G(t, x-y)  \psi(x) dx  \right) \varphi(y) dy  \\
     &=& \int_{\mathbb R^N} (e^{t \Delta}\psi)(y)  \varphi(y) dy.
  \end{eqnarray*}
  This proves the desired identity.
  \item [(ii)] $\int_{\mathbb R^N} \varphi (e^{t \Delta}\psi) \to  \int_{\mathbb R^N} \varphi \psi$ as $t \to 0$. In fact, since $\psi\in C_c^\infty(\mathbb R^N)$, we have 
      \begin{equation*}
       \|e^{t\Delta}\psi-\psi\|_{L^\infty(\mathbb R^N)} \rightarrow0, \quad\text{as }t \to 0. 
      \end{equation*}
       Consequently, 
       \begin{equation*} 
       (e^{t\Delta}\psi)(x)\rightarrow\psi(x), \quad\text{for every } x\in\mathbb R^N. 
       \end{equation*} 
       Hence, 
       \begin{equation*}
        \varphi(x)(e^{t\Delta}\psi)(x) \rightarrow \varphi(x)\psi(x), \quad\text{for every } x \in\mathbb R^N. 
        \end{equation*}
       It remains to find an integrable function, independent of $t$, which dominates $|\varphi(x)(e^{t\Delta}\psi)(x)|$. Suppose that $\operatorname{supp}(\psi)\subset B(0,R)=\{y \in \mathbb R^N;\quad |y|\leq R\}$, for some $R>0$. For $|x|\geq 2R$ and $|y|\leq R$, we have $|x-y|\geq |x|-|y|\geq\frac{|x|}{2}$. Therefore, 
       \begin{eqnarray}\label{eqn45}
       \left|(e^{t\Delta}\psi)(x)\right| &\leq& (4\pi t)^{-\frac N2} \int_{B(0,R)} e^{-\frac{|x-y|^2}{4t}}|\psi(y)| dy \nonumber\\
        &\leq& C_Rt^{-\frac N2} e^{-\frac{|x|^2}{16t}}, \quad \forall |x|\geq 2 R. 
        \end{eqnarray} 
       Moreover, for $|x|\geq\sqrt{8N}$, the function  $t\mapsto t^{-\frac N2}e^{-\frac{|x|^2}{16t}}$  is increasing on $(0,1]$. Hence, 
       \begin{equation}\label{eq451}
        t^{-\frac N2}e^{-\frac{|x|^2}{16t}} \leq e^{-\frac{|x|^2}{16}}, \quad \forall  0<t\leq1, ~  |x|\geq\sqrt{8N}.
       \end{equation}
       Set  $M=\max\{2R,\sqrt{8N}\}$.  Combining \eqref{eqn45} and \eqref{eq451}, we obtain, for $0<t\leq1$ and $|x|\geq M$, 
       \begin{equation*} 
       |\varphi(x)(e^{t\Delta}\psi)(x)| \leq C_R\|\varphi\|_{L^\infty_{\frac l\alpha}} |x|^{-\frac l\alpha} e^{-\frac{|x|^2}{16}}. 
       \end{equation*}
       On the other hand, for $|x| <M$, we use \begin{equation*}
         |(e^{t\Delta}\psi)(x)| \leq \|\psi\|_{L^\infty}, 
       \end{equation*}
          and consequently, 
       \begin{equation*} 
          |\varphi(x)(e^{t\Delta}\psi)(x)| \leq \|\varphi\|_{L^\infty_{\frac l\alpha}} \|\psi\|_{L^\infty} |x|^{-\frac l\alpha}. 
          \end{equation*}
       The dominated convergence theorem now proves the desired convergence.
\end{itemize}
Combining (i) and (ii), we conclude that $e^{t \Delta}\varphi \to \varphi$ as $t \to 0$  in $\mathcal D'(\mathbb R^N)$. Hence by \cite[Proposition 2.1 (i), p. 1110]{CDW-DCDS},  $e^{t \Delta}\varphi \to \varphi$ weak-$^*$ as $t \to 0$ in $\mathcal{B}_{\max\{1, C\}\|\varphi\|_{L_{\frac{l}{\alpha}}^\infty}}^*$, where $C$ is the constant obtained in \eqref{eq below} below and $\mathcal{B}_{\max\{1, C\}\|\varphi\|_{L_{\frac{l}{\alpha}}^\infty}}^*$ is defined in \eqref{BK}. Therefore, by \cite[Proposition 3.13 (iii), p. 63]{brezis},  $\|\varphi\|_{L_{\frac{l}{\alpha}}^\infty} \le \liminf_{t \to 0} \|e^{t \Delta}\varphi\|_{L_{\frac{l}{\alpha}}^\infty}$.
\end{proof}
\begin{rem}
{\rm It follows from \cite[Lemma 2.1]{CIT} that  
\begin{equation}\label{eq below}
  \|e^{t \Delta}\varphi\|_{L_{\frac{l}{\alpha}}^\infty} \le C \|\varphi\|_{L_{\frac{l}{\alpha}}^\infty}, \quad \forall t>0.
\end{equation}
Hence,  $\limsup_{t \to 0} \|e^{t \Delta}\varphi\|_{L_{\frac{l}{\alpha}}^\infty} \le C \|\varphi\|_{L_{\frac{l}{\alpha}}^\infty}$. Moreover, the argument used in the proof of Theorem \ref{thmupper1bis} yields  $\|\varphi\|_{L_{\frac{l}{\alpha}}^\infty} \le \liminf_{t \to 0} \|e^{t \Delta}\varphi\|_{L_{\frac{l}{\alpha}}^\infty}$. However, these two inequalities do not suffice to establish that $\|e^{t \Delta}\varphi\|_{L_{\frac{l}{\alpha}}^\infty} \to \|\varphi\|_{L_{\frac{l}{\alpha}}^\infty}$ as $t \to 0$. Consequently, the method employed in \cite[Proof of Theorem 1.6, p. 595]{TW4} does not apply directly.}
\end{rem}

\begin{proof}[Proof of Theorem \ref{thmupper1}]
By inequality \eqref{neccdbisbisbis}, we have
\begin{equation}\label{limsup}
  \lambda^{\frac{2 \alpha}{2+l}} T_{\max}(\lambda \varphi)  \leq  \left(\frac{2+l}{2 \alpha}\right)^{\frac{2}{2+ l}} \left[\frac{C_l}{\left({e}^{T_{\max}(\lambda \varphi)\Delta}\varphi \right)(0)}\right]^{\frac{2 \alpha}{2+ l}}.
\end{equation}
Passing to the $\limsup$ as $\lambda \to \infty$ in \eqref{limsup}, and using \eqref{limtmax} together with the fact that $e^{t\Delta}$ is a $C_0$-semigroup on $C_0(\mathbb{R}^N)$, we obtain the desired estimate \eqref{ineq45}.
\end{proof}
\begin{proof}[Proof of Theorem \ref{thmupper2bis}]
Let
\begin{equation}\label{phitilde}
\widetilde{\varphi}(x) =
\begin{cases}
\omega(x)\,|x|^{-\gamma}, & |x| \leq \varepsilon,\\
0, & |x| > \varepsilon.
\end{cases}
\end{equation}
Since $0\le \widetilde{\varphi}\le \varphi$, we have
$T_{\max}(\lambda\varphi)\le T_{\max}(\lambda\widetilde{\varphi})$.
Thus, it suffices to estimate $T_{\max}(\lambda\widetilde{\varphi})$.
Moreover, as in \cite[p. 595]{TW4}, We have
\begin{equation}
\label{gammascale}
\|{\rm e}^{t\Delta}\tilde\varphi\|_\infty = \|D_{\sqrt t}{\rm e}^{t\Delta}\tilde\varphi\|_\infty
=\|{\rm e}^{\Delta}D_{\sqrt t}\tilde\varphi\|_\infty
= t^{-\frac{\gamma}{2}}\|{\rm e}^{\Delta}[t^{\frac{\gamma}{2}}D_{\sqrt t}\tilde\varphi]\|_\infty.
\end{equation}
Since $t^{\frac{\gamma}{2}}D_{\sqrt t}\tilde\varphi \to \omega|\cdot|^{-\gamma}$ in ${\mathcal D'}(\R^N)$
as $t \to 0$, it follows by \cite[Proposition 3.8 (i), page 1123]{CDW-DCDS}, that
\begin{equation*}
  t^{\frac{\gamma}{2}}  \left\|e^{t \Delta}\widetilde{\varphi} \right\|_{L^\infty} \to \left\| e^{\Delta}(\omega |\cdot|^{-\gamma})\right\|_{L^\infty}, \quad \mbox{as} \quad t\rightarrow 0.
\end{equation*}
Since $T_{\max}(\lambda\varphi)\to 0$ as $\lambda\to\infty$ by \eqref{limtmax}, we deduce from \eqref{neccdinftyinftybis} that
\begin{equation*}
  \limsup_{\lambda \to \infty}\lambda^{\frac{2 \alpha}{2+l}} \left(T_{\max}(\lambda \widetilde{\varphi}) \right)^{1- \frac{\alpha \gamma}{2+ l}}  \leq  \left(\frac{2+l}{2 \alpha}\right)^{\frac{2}{2+ l}} \left[\frac{C_l}{\left\|e^{\Delta}(\omega |\cdot|^{-\gamma}) \right\|_{L^\infty}}\right]^{\frac{2 \alpha}{2+ l}},
\end{equation*}
which completes the proof.
\end{proof}
\begin{proof}[Proof of Theorem \ref{thmupper2}]
It suffices to estimate $T_{\max}(\lambda \widetilde{\varphi})$ from above, where $\widetilde{\varphi}$ is defined in \eqref{phitilde}. For all $t>0$, we have
\begin{equation}\label{eqdfr}
   t^{\frac{\gamma}{2}}  \left(e^{t \Delta}\widetilde{\varphi} \right)(0) =  \left[e^{ \Delta} \left(t^{\frac{\gamma}{2}} D_{\sqrt t}\widetilde{\varphi}\right)\right](0).
\end{equation}
Since  $ t^{\frac{\gamma}{2}} D_{\sqrt t}\widetilde{\varphi}\to \omega |\cdot|^{-\gamma}$ in $\mathcal D'(\mathbb R^N)$ as $t \to 0$, it follows by \cite[Proposition 3.8 (i), p. 1123]{CDW-DCDS}, that 
\begin{equation*}
  e^\Delta \left( t^{\frac{\gamma}{2}} D_{\sqrt{t}} \widetilde{\varphi} \right) \to e^{\Delta}(\omega |\cdot|^{-\gamma}) \quad \text{in } C_0(\mathbb{R}^N), \quad \mbox{as} \quad t\rightarrow 0,
\end{equation*}
and consequently,
\begin{equation}\label{ineq4578}
  t^{\frac{\gamma}{2}}  \left(e^{t \Delta}\widetilde{\varphi} \right)(0) \to e^{\Delta}(\omega |\cdot|^{-\gamma})(0), \quad \mbox{as} \quad t\rightarrow 0.
\end{equation}
Finally, since $T_{\max}(\lambda \varphi)\rightarrow 0$ as $\lambda\rightarrow \infty$, combining \eqref{neccdbisbisbis} with \eqref{ineq4578} yields
\begin{equation*}
  \limsup_{\lambda \to \infty}\lambda^{\frac{2 \alpha}{2+l}} \left(T_{\max}(\lambda \widetilde{\varphi}) \right)^{1- \frac{\alpha \gamma}{2+ l}}  \leq  \left(\frac{2+l}{2 \alpha}\right)^{\frac{2}{2+ l}} \left[\frac{C_l}{e^{\Delta}(\omega |\cdot|^{-\gamma})(0)}\right]^{\frac{2 \alpha}{2+ l}},
\end{equation*}
which concludes the proof.
\end{proof}

\begin{proof}[Proof of Theorem \ref{thmupper3bis}]
Let
\begin{equation}\label{phi2tilde}
  \widetilde{\widetilde{\varphi}}(x) =
\begin{cases}
0, & |x| < R, \\
\omega(x)\,|x|^{-\gamma}, & |x| \ge R.
\end{cases}
\end{equation}
Since $\varphi \geq \widetilde{\widetilde{\varphi}}$, it is sufficient to derive an upper bound for $T_{\max}(\lambda \widetilde{\widetilde{\varphi}})$.   As in \cite[p. 596]{TW4}, the calculation in \eqref{gammascale} gives
\begin{equation*}
%\label{gammascale2}
\|{\rm e}^{t\Delta}\tilde{\tilde\varphi}\|_{L^\infty}
= t^{-\frac{\gamma}{2}}\|{\rm e}^{\Delta}[t^{\frac{\gamma}{2}}D_{\sqrt t}\tilde{\tilde\varphi}]\|_{L^\infty}.
\end{equation*}
Moreover, $t^{\frac{\gamma}{2}}D_{\sqrt t}\tilde{\tilde\varphi} \to \omega|\cdot|^{-\gamma}$
as $t \to \infty$ in ${\mathcal D'}(\R^N).$ It follows, by \cite[Proposition 3.8 (i), page 1123]{CDW-DCDS}  that
\begin{equation}
\label{equal1}
t^{\frac{\gamma}{2}}\|{\rm e}^{t\Delta}\tilde{\tilde\varphi}\|_{L^\infty}
= \|{\rm e}^{\Delta}[t^{\frac{\gamma}{2}}D_{\sqrt t}\tilde{\tilde\varphi}]\|_{L^\infty}
\to \|{\rm e}^{\Delta}[\omega|\cdot|^{-\gamma}]\|_{L^\infty}, \quad \mbox{as}\quad t \to \infty.
\end{equation}
If $0 < t < T_{\max}(\lambda\tilde{\tilde\varphi})$, then by \eqref{neccdinfty} we must have
\begin{equation}
\label{gammascale3}
\frac{2}{2+l} \alpha \lambda^\alpha t^{\frac{2+l}{2}- \frac{\alpha\gamma}{2}} [t^{\frac{\gamma}{2}}\|e^{t\Delta}\tilde{\tilde\varphi}\|_{L^\infty} ]^\alpha \leq C_l^\alpha.
\end{equation}
It follows that if $\gamma < \frac{2+l}{\alpha}$, then $T_{\max}(\lambda\tilde{\tilde\varphi}) < \infty$
for all $\lambda > 0$. This is of course a consequence of Fujita's result (including the limiting case)
if $\alpha \le \frac{2+l}{N}$. (See for example \cite{Qi}.) If $\alpha > \frac{2+l}{N}$,  this is a consequence of Theorem \ref{thmblowupHénon} and Remark \ref{compact} (iv). On the other hand, it follows from Theorem \ref{thm lower bound 2} when $\frac{l}{\alpha} < \gamma < \min \{N, \frac{2+l}{\alpha}\}$,  and from Theorem \ref{thm lower bound 2bis} when $\gamma=\frac{l}{\alpha}$,  that $T_{\max}(\lambda \widetilde{\widetilde{\varphi}}) \to \infty$ as $\lambda \to 0.$
In these cases, since $T_{\max}(\lambda\tilde{\tilde\varphi})\to\infty$ as $\lambda\to 0,$ putting $t = T_{\max}(\lambda\tilde{\tilde\varphi})$  in \eqref{gammascale3}
and letting $\lambda \to 0$ and using \eqref{equal1}, we obtain
\begin{equation*}
\limsup_{\lambda \to 0} \lambda^{(\frac{2+l}{2\alpha} - \frac{\gamma}{2})^{-1}} T_{\max}(\lambda\tilde{\tilde\varphi})
\le\left[\frac{C_l}{\left(2\alpha/(2+l)\right)^{1/\alpha}\left\|e^{\Delta}(\omega |\cdot|^{-\gamma})\right\|_{L^\infty}}\right]^{(\frac{2+l}{2\alpha} - \frac{\gamma}{2})^{-1}}.
\end{equation*}
This completes the proof.
\end{proof}

\begin{proof}[Proof of Theorem \ref{thmupper3}]
Let $\widetilde{\widetilde{\varphi}}$ be as in \eqref{phi2tilde}. It is sufficient to derive an upper bound for $T_{\max}(\lambda \widetilde{\widetilde{\varphi}})$. First, observe that if $\gamma < \frac{2+l}{\alpha}$, then $T_{\max}(\lambda \widetilde{\widetilde{\varphi}}) < \infty$ for all $\lambda > 0$, as a direct consequence of Theorem \ref{thmblowup}, and Remark \ref{compact} (iv). Next, replacing $\widetilde{\varphi}$ by $\widetilde{\widetilde{\varphi}}$ in \eqref{eqdfr}, we obtain
\begin{equation*}
t^{\frac{\gamma}{2}} \left(e^{t \Delta}\widetilde{\widetilde{\varphi}}\right)(0) = \left[e^{\Delta} \left(t^{\frac{\gamma}{2}} D_{\sqrt{t}} \widetilde{\widetilde{\varphi}}\right)\right](0).
\end{equation*}
Moreover, as $t \to \infty$, we have $t^{\frac{\gamma}{2}} D_{\sqrt{t}} \widetilde{\widetilde{\varphi}} \to \omega |\cdot|^{-\gamma}$ in $\mathcal{D}'(\mathbb{R}^N).$ Therefore, by \cite[Proposition 3.8, p. 1123]{CDW-DCDS}, it follows that
\begin{equation*}
t^{\frac{\gamma}{2}} \left(e^{t \Delta}\widetilde{\widetilde{\varphi}}\right)(0) \to \left(e^{\Delta}(\omega |\cdot|^{-\gamma})\right)(0), \quad \text{as } t \to \infty.
\end{equation*}
It follows from Theorem \ref{thm lower bound 2} that  $T_{\max}(\lambda \widetilde{\widetilde{\varphi}}) \to \infty$ as $\lambda \to 0$. Combining this with \eqref{neccdbisbisbis}, implies that
\begin{equation*}
  \limsup_{\lambda \to 0}\lambda^{\frac{2 \alpha}{2+l}} \left(T_{\max}(\lambda \widetilde{\widetilde{\varphi}}) \right)^{1- \frac{ \alpha \gamma}{2+ l}}  \leq  \left(\frac{2+l}{2 \alpha}\right)^{\frac{2}{2+ l}} \left[\frac{C_l}{e^{\Delta}(\omega |\cdot|^{-\gamma})(0)}\right]^{\frac{2 \alpha}{2+ l}}.
\end{equation*}
\end{proof}


\begin{thebibliography}{99}
\bibitem{AF}{R. A.  Adams and J. J. Fournier, {\it Sobolev spaces}, Elsevier, 2003.}

\bibitem{BL}{C. Bandle and H. A. Levine, {\it On the existence and nonexistence of global solutions of reaction-diffusion equations in sectorial domains}, Trans. AMS, 316 (1989), 595--622.}
    
\bibitem{BenSlimene}{B. Ben Slimene, {\it Asymptotically self-similar global solutions for Hardy-H\'enon parabolic systems},  Differ. Equ. Appl., 11 (2019), 439--462.}
    
\bibitem{BTW} {B.   Ben Slimene,  S.  Tayachi and  F. B.  Weissler, {\it  Well-posedness, global existence and large time behavior for Hardy-H\'enon parabolic equations},  Nonlinear Anal.,  152 (2017),  116--148.}
    
\bibitem{brezis}{H. Br\'ezis, {\it Functional analysis, Sobolev spaces and partial differential equations}, Universitext, Springer, New York, 2011.}

\bibitem{CFL}{ R. Castillo, R. Freire and M. Loayza, {\it Global existence versus blow-up for a Hardy-H\'enon parabolic equation on arbitrary domains}, J. Differential Equations, 429 (2025), 427--459.}

\bibitem{CDW-DCDS}{T. Cazenave, F. Dickstein and F. B. Weissler, {\it Universal solutions of the heat equation in $\Rd$,}  Discrete Contin. Dynam. Systems, 9 (2003), 1105--1132.}

\bibitem{Ch}{N. Chikami, {\it Composition estimates and well-posedness for Hardy-Hénon parabolic equations in Besov spaces}, J. Elliptic Parabol. Equ., 5 (2019) 215--250.}

\bibitem{CIT2}{N. Chikami, M. Ikeda and K. Taniguchi, {\it Well-posedness and global dynamics for the critical Hardy-Sobolev parabolic equation}, Nonlinearity, 34 (2021) 8094--8142.}

\bibitem{CIT}{N. Chikami, M. Ikeda and K. Taniguchi, {\it Optimal well-posedness and forward self-similar solution for the Hardy-H\'enon parabolic equation in critical weighted Lebesgue spaces}, Nonlinear Analysis,  222 (2022), 112931.}

\bibitem{CITT}{N. Chikami, M. Ikeda, K. Taniguchi and S. Tayachi, {\it Asymptotically self-similar global solutions for Hardy-Hénon parabolic equations}, J. Differential Equations, 466 (2026), Paper No. 114326, 70 pp.}

\bibitem{CITT2}{N. Chikami, M. Ikeda, K. Taniguchi and  S. Tayachi, {\it Unconditional uniqueness and non-uniqueness for Hardy-Hénon parabolic equations}, 
Math. Ann., 390 (2024), 3765–3825.}

\bibitem{FT}{S. Filippas and A. Tertikas, {\it On similarity solutions of a heat equation with a nonhomogeneous nonlinearity}, J. Differential Equations, 165 (2000),  468--492.}

\bibitem{Folland}{G. B. Folland, {\it Real analysis: modern techniques and their applications},  volume 40. John Wiley \& Sons, 1999.}

\bibitem{F}{H. Fujita, {\it On the blowing up of solutions of the Cauchy problem for $u_{t}=\Delta u+u^{1+\alpha}$},  J. Fac. Sci. Univ. Tokyo Sect. I,  13  (1966), 109--124.}

\bibitem{Fujita2}{ H. Fujita, {\it On some nonexistence and nonuniqueness theorems for nonlinear parabolic equations}, Proc. Sympos. Pure Math., 18 (1968), 138--161.}

\bibitem{GW}{C. Gui and X. Wang, {\it Life span of solutions of the Cauchy problem for a semilinear heat equation}, J. Differential Equations,  115 (1995), 166--172.}

\bibitem{GLS1}{J.-S. Guo, C.-S. Lin and M. Shimojo, {\it Blow-up behavior for a parabolic equation with spatially dependent coefficient}, Dynam. Systems Appl., 19 (2010),  415--433.}

\bibitem{GLS2}{J.-S. Guo, C.-S. Lin and M. Shimojo, {\it  Blow-up for a reaction-diffusion equation with variable coefficient}, Appl. Math. Lett., 26
8 (2013), 150--153.}

\bibitem{GS}{J.-S. Guo and M. Shimojo, {\it Blowing up at zero points of potential for an initial boundary value problem}, Commun. Pure Appl. Anal., 10 (2011), 161--177.}

\bibitem{GSouplet}{J.-S. Guo and Ph. Souplet, {\it Excluding blowup at zero points of the potential by means of Liouville-type theorems}, J. Differential Equations, 265 (2018), 4942--4964.}

\bibitem{HisaIshige}{K. Hisa and K. Ishige, {\it Existence of solutions for a fractional semilinear parabolic equation with singular initial data}, Nonlinear Anal., 175 (2018), 108--132.}


\bibitem{IL}{R. G. Iagar and  Ph. Lauren\c cot, {\it Sharp non-existence threshold for a parabolic Hardy-H{é}non equation with quasilinear diffusion}, 	arXiv:2508.06164 (math.AP).}

\bibitem{LeeNi}{T. Y. Lee and  W. M. Ni, {\it Global existence, large time behavior and life span of solutions of a semilinear parabolic Cauchy problem},
Trans. Amer. Math. Soc., 333 (1992), 365--378.}

\bibitem{LM}{H. A. Levine and P. Meier, {\it The value of the critical exponent for reaction-diffusion equations in cones}, Arch. Ration. Mech. Anal.,
 109 (1990), 73--80.}

\bibitem{MY}{N. Mizoguchi and E. Yanagida, {\it Blowup and life span of solutions for a semilinear parabolic equation}, SIAM J. Math. Anal., 29 (1998), 1434--1446.}

\bibitem{MTW}{H. Mouajria, S. Tayachi and F. B. Weissler, {\it The heat equation on sectorial domains, highly singular initial values and applications}, J. Evol. Equ., 16 (2016), 341--364.}

\bibitem{Pinsky}{R. G.  Pinsky, {\it  Existence and Nonexistence of Global Solutions for  $u_t= \Delta u+ a(x) u^p$ in $\mathbb R^d$}, J. Differential
  Equations,  133 (1997), 152--177.}

\bibitem{Qi}{Y. W.  Qi, {\it  The critical exponents of parabolic equations and blow-up in  $\mathbb R^n$}, Proc. Roy. Soc. Edinburgh Sect.,  128A (1998), 123--136.}

\bibitem{R}{P. Rouchon, {\it Blow-up of solutions of nonlinear heat equations in unbounded domains for slowly decaying initial data}, ZAMP, 52 (2001), 1017--1032.}

\bibitem{SW}{ Ph. Souplet and F. B. Weissler, {\it Self-similar subsolutions and blow-up for nonlinear parabolic equations}, J. Math. Anal. Appl., 212 (1997), 60--74.}

\bibitem{Tayachi}{S. Tayachi, {\it Uniqueness and Non-uniqueness  of solutions for critical Hardy-H\'enon parabolic equations}, J. Math. Anal. Appl. 488 (2020), 123976.}
    
\bibitem{Tayachi1}{  S.  Tayachi, {\it Life-span of solutions for a nonlinear parabolic system},  NoDEA Nonlinear Differ. Equ. Appl., 31 (2024),  67--111.}

\bibitem{TW4}{S. Tayachi and  F. B. Weissler, {\it New life-span results for the nonlinear heat equation}, J. Differ. Equ., 373 (2023) 564--625.}

\bibitem{TW5}{S. Tayachi and  F. B. Weissler, {\it The nonlinear heat equation involving highly singular initial values and new blowup and life span results}, J. Elliptic Parabol. Equ., 4 (2018), 141--176.}

\bibitem{Tsusaka}{Y. Tsusaka, {\it Local existence and nonexistence of solutions to the Hardy parabolic equation with general nonlinearity}, arXiv:2601.22520 (math.AP).}


\bibitem{Wang}{ X.  Wang, {\it On the Cauchy problem for reaction-diffusion equations}, Trans. Amer. Math. Soc., 337 (1993), 549--590.}

\bibitem{W4}{F. B. Weissler, { \it $L^p-$Energy and blow-up for a semilinear heat equation}, Proceedings of Symposia in Pure Mathematics, 45 (1986), Part 2, 545--551}.

\end{thebibliography}
\end{document}